\documentclass[12pt, reqno, a4paper,oneside]{amsart}

\usepackage{graphicx, epsfig, psfrag}
\usepackage{amsfonts,amsmath,amssymb,amsbsy,amsthm}
\usepackage{bm}
\usepackage{color}
\usepackage{float}
\usepackage{hyperref}
\usepackage{mathrsfs}
\usepackage{longtable}
\usepackage[top=3cm,bottom=3cm,left=3cm,right=3cm]{geometry}
\usepackage{epstopdf}
\usepackage{mathabx}
\usepackage{stmaryrd}
\usepackage{ulem}
\usepackage{scalerel}
\usepackage{enumitem}
\usepackage{subfigure}

\usepackage{environments}
\numberwithin{theorem}{section}
\numberwithin{equation}{section}

\newtheorem{algorithm}{Algorithm}[section]

\renewcommand{\cases}[1]{\left\{ \begin{array}{rl} #1 \end{array} \right.}

\newcommand{\mymat}[1]{\left[ \begin{matrix} #1 \end{matrix} \right]}

\def\XXint#1#2#3{{\setbox0=\hbox{$#1{#2#3}{\int}$ }
\vcenter{\hbox{$#2#3$ }}\kern-.6\wd0}}

\def\b{\big}
\def\B{\Big}

\def\sep{\,|\,}
\def\bsep{\,\b|\,}

\def\R{\mathbb{R}}

\def\Z{\mathbb{Z}}

\def\dx{\,{\textrm{d}}x}

\def\dt{\,{\textrm{d}}t}
\def\ds{\,{\textrm{d}}s}

\def\pp{\partial}

\def\<{\langle}
\def\>{\rangle}

\def\mA{{\textsf{A}}}
\def\mB{\textsf{B}}

\def\mF{\textsf{F}}

\def\mJ{\textsf{J}}
\def\mP{\textsf{P}}

\def\mR{\textsf{R}}

\newcommand{\Dc}[1]{\D_{#1}}
\def\D{\nabla}
\def\del{\delta}
\def\ddel{\delta^2}

\def\c{\textrm{c}}
\def\h{\textrm{h}}

\def\a{\textrm{a}}
\def\c{\textrm{c}}
\def\ac{\textrm{ ac}}
\def\i{\textrm{i}}

\def\L{\Lambda}

\def\Us{\mathscr{U}}

\def\Ush{\dot{\Us}^{1,2}}

\def\Usc{\Us^c}

\def\E{\mathscr{E}}
\def\Ea{\E^\a}
\def\Eb{\E^\textrm{b}}

\def\Ei{\E^\i}
\def\Ec{\E^\c}
\def\Eh{\E^\h}

\def\L{\Lambda}
\def\La{\L^{\a}}
\def\Li{\L^{\i}}

\def\Nhd{\mathcal{N}}
\def\rcut{r_{\textrm{cut}}}
\def\Rg{\mathscr{R}}
\def\Rgp{\Rg^{+}}
\def\vsig{\varsigma}
\def\T{\mathcal{T}}

\def\T{\mathcal{T}}
\def\Th{\mathcal{T}_h}

\def\Fh{\mathscr{F}_h}
\def\UsT{\Us_h}

\def\Ia{I_\a}

\def\vor\textrm{vor}
\def\s{\sigma}
\def\sh{\sigma^h}
\def\sa{\sigma^\a}
\def\sc{\sigma^\c}

\def\PO{\textrm{P}_0}

\definecolor{lzcol}{rgb}{0.7, 0, 0}
\definecolor{mlcol}{rgb}{0, 0.7, 0}
\definecolor{todocol}{rgb}{0.0, 0.4, 0.0}

\begin{document}
\title[A Posteriori Estimate and Adaptivity for Finite Range A/C 2D]{A Posteriori Error Estimate and Adaptive Mesh Refinement Algorithm for Atomistic/Continuum Coupling with Finite Range Interactions in Two Dimensions}

\author{M. Liao}
\address{Mingjie Liao\\
Department of Applied Mathematics and Mechanics\\
University of Science and Technology Beijing\\
No. 30 Xueyuan Road, Haidian District\\
Beijing 100083\\
China}
\email{mliao@xs.ustb.edu.cn}

\author{P. Lin}
\address{Ping Lin\\
Department of Mathematics\\
University of Dundee\\
Dundee, DD1 4HN, Scotland\\
United Kingdom}
\email{plin@maths.dundee.ac.uk}

\author{L. Zhang}
\address{Lei Zhang \\ School of Mathematical Sciences,
  Institute of Natural Sciences, and Ministry of Education Key
  Laboratory of Scientific and Engineering Computing (MOE-LSC) \\
  Shanghai Jiao Tong University \\ 800 Dongchuan Road \\ Shanghai
  200240 \\ China}
\email{lzhang2012@sjtu.edu.cn}

\thanks{ML and PL were partially supported by National Natural Science Foundation of China grant 91430106. LZ was partially supported by National Natural Science Foundation of China grant 11471214, 11571314 and the One Thousand Plan of China for young scientists.}

\begin{abstract}
In this paper, we develop the residual based a posteriori error estimates and the corresponding adaptive mesh refinement algorithm for atomistic/continuum (a/c) coupling with finite range interactions in two dimensions. We have systematically derived a new explicitly computable stress tensor formula for finite range interactions. In particular, we use the geometric reconstruction based consistent atomistic/continuum (GRAC) coupling scheme, which is optimal if the continuum model is discretized by $P^1$ finite elements. The numerical results of the adaptive mesh refinement algorithm is consistent with the optimal a priori error estimates.

\end{abstract}


\maketitle

\section{Introduction}
\label{sec:introduction}
Atomistic/continuum (a/c) coupling methods are a class of computational multiscale methods for crystalline solids with defects that aim to optimally balance the accuracy of the atomistic model and the efficiency of the continuum model \cite{Ortiz:1995a, Shenoy:1999a, Gumbsch:1989}. The construction and analysis of different a/c coupling methods have attracted considerable attention in the research community in recent years. Rigorous a prior analysis and systematic benchmark has been done in, for example, \cite{LinP:2006a, OrtnerShapeev:2011, LiOrShVK:2014, COLZ2013, OrZh:2016, MiLu:2011}. We refer readers to \cite{Miller:2008, LuOr:acta} for a review. The study of a/c coupling methods has not only provided an analytical framework for the prototypical problems\cite{Ehrlacher:2016}, but also opened avenue for coupling schemes in more complicated physical situations \cite{QMMM:2016,QMMM:2017, Fang:2018}.

Like many multiscale methods dealing with defects or singularities, adaptivity is the key for the efficient implementation of a/c coupling methods. In contrast to the a priori analysis, the development of a posteriori analysis for a/c coupling methods are still lagging behind. Although heuristic methods have been proposed in the engineering literature \cite{Kochman:2016, prud06,Shenoy:1999a}. Previous mathematical justifications are largely limited to one dimension cases \cite{prud06, arndtluskin07b, arndtluskin07c}. In particular, the residual based a posteriori error bounds for a/c coupling schemes are first derived in \cite{OrtnerSuli:2008a, Ortner:qnl.1d, OrtnerWang:2014} by Ortner et al. in one dimension. 

In \cite{APEst:2017}, we carried out a rigorous a posteriori analysis of the residual, the stability constant, and the error bound, for a consistent atomistic/continuum coupling method \cite{PRE-ac.2dcorners} with nearest neighbour interactions in two dimensions. Corresponding adaptive mesh refinement algorithm was designed and implemented based on the a posteriori error estimates, and the convergence rate with respect to degrees of freedom is the same as optimal a priori error estimates. This is the first rigorous a posteriori analysis for a/c coupling method in two dimensions. With the a posteriori error estimates and the adaptive algorithm, we can not only automatically move the a/c interface and adjust the discretization of the continuum region, but also change the size of the computational domain. We have also introduced the so-called ``stress tensor correction" technique, which distinguish the essential difference of high dimensional results compared with previous one dimensional results. 

In this paper, we treat the more general case of a/c coupling with finite range interactions, which is physically more relevant and computationally more involved. The a priori analysis of GRAC scheme has been extended from nearest neighbour case in \cite{PRE-ac.2dcorners} to finite range interactions in \cite{COLZ2013}. $\ell^{1}$-minimization is introduced to resolve the issue of non-uniqueness of reconstruction parameters, and a stabilisation mechanism is proposed in \cite{2013-stab.ac} to reduce the stability gap between the a/c coupling scheme and the original atomistic model.

The analytical framework for both a priori analysis and a posteriori analysis of a/c coupling methods strongly relies on the stress based formulation. In \cite{Or:2011a,OrtnerTheil2012}, an explicit formulation of stress tensor is proposed based on a mollified version of line measure supported on the interaction bonds, thence one can obtain an integral representation of finite differences to further derive an integral representation of the first variation of the interaction energy. This representation greatly simplifies the expression of the stress tensor and plays a significant role in the a priori analysis. However, the obtained stress tensor is a function of the continuous space variable, therefore it is difficult to compute in practice, and not suitable for the a posteriori estimates and adaptive computation. 

In this paper, we derive a novel expression of the stress tensor for finite range interactions, which is new to the best of our knowledge. The stress tensor is piecewise constant and only depends on a local neighbourhood, therefore it is computable and the assembly cost is linear with respect to the number of bonds. This stress formulation allows for the convenient derivation of a posterior error estimates and efficient implementation of the adaptive algorithms.

The paper is organized as follows. We set up the atomistic to continuum (a/c) coupling models for point defects in \S~\ref{sec:formulation}, and introduce the general GRAC formulation in \S~\ref{sec:grac}. In \S~\ref{sec:stress}, we present the stress formulation and the stress tensor assembly algorithm for finite range interactions. In \S~\ref{sec:numerics} the rigorous a posteriori error estimates based adaptive algorithm is deduced and complemented by numerical experiments. We draw conclusions and make suggestions for future research in \S~\ref{sec:conclusion}. Some auxiliary results are given in the Appendix \S~\ref{sec:appendix}.


%

\def\Rdef{R^{\textrm{def}}}
\def\Rg{\mathcal{R}}
\def\Rgnn{\mathcal{N}}
\def\rcut{r_{\textrm{cut}}}
\def\Lhom{\L^{\textrm{hom}}}
\def\Ddef{D^{\textrm{def}}}
\def\Ldef{\L^{\textrm{def}}}

\def\Nh{\mathcal{N}_h}
\def\Ush{\Us_h}
\def\Ra{R^\a}
\def\Rb{R^{\textrm{b}}}
\def\Rc{R^\c}
\def\Eb{\E^{\textrm{b}}}
\def\dof{{\textrm{DOF}}}
\def\Omh{\Omega_h}
\def\Thr{{\T_{h,R}}}
\def\vor{\textrm{vor}}
\def\Uhr{\Us_{h,R}}

\def\Ta{\T_\a}
\def\Th{\T_{\textrm{h}}}
\def\sh{\sigma^{\textrm{h}}}

\section{Model Formulation}
\label{sec:formulation}
In this section, We setup an atomistic model for crystal defects in an infinite lattice in the spirit of \cite{Ehrlacher:2016} in \S~\ref{sec:formulation:atm} and then introduce the Cauchy-Born continuum model in \S~\ref{sec:formulation:continuum}. We give a generic form of a/c coupling schemes in \S~\ref{sec:formulation:ac}. 

\subsection{Atomistic Model}
\label{sec:formulation:atm}
\def\Rdef{R^{\textrm{def}}}
\def\Rg{\mathcal{R}}
\def\Rgnn{\mathcal{N}}
\def\rcut{r_{\textrm{cut}}}
\def\Lhom{\L^{\textrm{hom}}}
\def\Ddef{D^{\textrm{def}}}
\def\Ldef{\L^{\textrm{def}}}
\subsubsection{Atomistic lattice and defects}
\label{sec:formulation:atm:lattice}

Given $d \in \{2, 3\}$, $\mA \in \R^{d \times d}$ non-singular, $\Lhom := \mA \Z^d$ is the \textit{homogeneous reference lattice} which represents a perfect \textit{single lattice crystal} formed by identical atoms. $\L\subset \R^d$ is the \textit{reference lattice} with some local \textit{defects}. The mismatch between $\L$ and $\Lhom$ represents possible defects $\Ldef$, which are contained in some localized \textit{defect cores} $\Ddef$ such that the atoms in $\L\setminus \Ddef$ do not interact with defects $\Ldef$. Vacancy, interstitial and impurity are different types of possible point defects.

\subsubsection{Lattice function and lattice function space}
\label{sec:formulation:atm:function}

Given $m\in\{1,2,3\}$, denote the set of vector-valued \textit{lattice functions} by 
\begin{displaymath}
\Us := \{v: \L\to \mathbb{R}^m \}.
\end{displaymath}

A \textit{deformed configuration} is a \textit{lattice function} $y \in \Us$. Let $x$ be the identity map, the \textit{displacement} $u\in \Us$ is defined by $u(\ell) = y(\ell)-x(\ell) = y(\ell) - \ell $ for any $\ell\in\L$. 

For each $\ell\in \L$, we prescribe an \textit{interaction neighbourhood}  $\Nhd_{\ell} := \{ \ell' \in \L \sep 0<|\ell'-\ell| \leq \rcut \}$ with some cut-off radius $\rcut$. The \textit{interaction range} $\Rg_\ell := \{\ell'-\ell \sep \ell'\in \Nhd_\ell\}$ is defined as the union of lattice vectors defined by the finite differences between lattice points in $\Nhd_{\ell}$ and $\ell$. Define the ``finite difference stencil'' $Dv(\ell):= \{D_\rho v(\ell)\}_{\rho \in \Rg_\ell} :=\{v(\ell+\rho)-v(\ell)\}_{\rho \in \Rg_\ell}$.

The homogeneous lattice $\Lhom = \mA\mathbb{Z}^d$ naturally induces a simplicial micro-triangulation $\T^{\a}$. In two dimensions, $\T^{\a} = \{\mA\xi + \hat{T}, \mA\xi-\hat{T}|\xi\in \mathbb{Z}^2\}$, where $\hat{T} = {\textrm{conv}}\{0, e_1, e_2\}$. Let $\bar{\zeta}\in W^{1, \infty}(\Lhom; \R)$ be the $P_1$ nodal basis function associated with the origin; namely, $\bar{\zeta}$ is piecewise linear with respect to $\T^\a$, and $\bar{\zeta}(0) = 1$ and $\bar{\zeta}(\xi)=0$ for $\xi\neq 0$ and $\xi\in \Lhom$. The nodal interpolant of $v\in \Us$ can be written as 
\begin{displaymath}
	\bar{v}(x):=\sum_{\xi\in\Z^d}v(\xi)\bar{\zeta}(x-\xi).
\end{displaymath}

\def\Use{\Us^{1,2}}

We can introduce the discrete homogeneous Sobolev spaces
\begin{displaymath}
	\Use :=\{u\in \Us|\nabla \bar{u}\in L^2\},
\end{displaymath}
with semi-norm $\|\nabla \bar{u}\|_{L^2}$.

\subsubsection{Interaction potential}
\label{sec:formulation:atm:potential}

We consider the general multibody interaction potential of the \textit{generic pair functional form} \cite{TadmorMiller:2012}, which includes the widely used potentials such as EAM (Embedded Atom Method) potential \cite{Daw1984a} and Finnis-Sinclair model \cite{Finnis1984}. Namely, the potential is a function of the distances between atoms within interaction range and has no angular dependence. For example, for each $\ell \in \L$, let $V_\ell(y)$ denote the \textit{site energy}
associated with the lattice site $\ell \in \L$, the EAM potential reads, 

\begin{align}
  V_{\ell}(y) := & \sum_{\ell' \in \Nhd_{\ell}} \Phi(|y(\ell)-y(\ell')|) + F\B(
  {\textstyle \sum_{\ell' \in \Nhd_{\ell}} \psi(|y(\ell)-y(\ell')|)} \B),\nonumber\\
    = &\sum_{\rho \in \Rg_{\ell}} \Phi\b(|D_\rho y(\ell)|\b) + F\B(
  {\textstyle \sum_{\rho \in \Rg_{\ell}}} \psi\b( |D_\rho y(\ell)|\b) \B).   
    \label{eq:eam_potential}
\end{align}
with the pair potential $\Phi$, the electron density function $\psi$ and the embedding function $F$. 


\begin{remark}
\label{rem:pairfuncational}
For convenience, with a slight abuse of notation, we will use $V_{\ell}(D_\rho y)$,  $V_{\ell}(|D_\rho y|)$ instead of $\widehat{V}_{\ell}(\{D_\rho y(\ell) \}_{\rho\in \Rg_\ell})$,  $\widetilde{V}_{\ell}(\{|D_\rho y (\ell) |\}_{\rho\in \Rg_\ell})$ when there is no confusion in the context.
\end{remark}

We assume that the potential $V_\ell(y)\in C^k((\R^d)^{\Rg_\ell}), k \geq 2$. We also assume that $V_\ell(y)$ is \textit{homogeneous} outside the defect region $\Ddef$, namely, $V_\ell = V$ and $\Rg_\ell = \Rg$ for $\ell \in \Lambda \setminus \Ddef$. Furthermore, $V$ and $\Rg$ have the following point symmetry: $\Rg = -\Rg$, and $V(\{-g_{-\rho}\}_{\rho\in\Rg}) = V(g)$. 

For an infinite lattice, assume the macroscopic applied strain is $\mB\in\R^{d\times d}$, we redefine the potential $V_\ell(y)$ as the difference $V_\ell(y) - V_\ell(y_\mB)$. We denote the energy functional $\E(y)$ as the infinite sum of the redefined potential over $\Lambda$, which is well-defined for $y-y^{\mB}\in\Use$ \cite{Ehrlacher:2016},  
\begin{equation}
\label{eqn:Ea}
  \E(y) = \sum_{\ell \in \L} V_\ell(y) 
\end{equation}

Under the above conditions, the goal of the atomistic problem is to find a \textit{strongly stable}
equilibrium $y$, such that, 
\begin{equation}
  \label{eq:min}
  y \in \arg\min \b\{ \E(y) \bsep y-y^B\in \Use \b\}.
\end{equation}
$y$ is \textit{strongly stable} if there exists $c_0 > 0$ such that
\begin{displaymath}
\< \ddel \E(y) v, v \> \geq c_0 \| \nabla v \|_{L^2}^2, \quad \forall v \in \Us^{1,2}.
\end{displaymath}

\subsection{Continuum model}
\label{sec:formulation:continuum}
From the atomistic model, a continuum model can be derived by coarse graining, and computationally it allows for the reduction of degrees of freedom when the deformation is smooth. A typical choice in the multi-scale context is the Cauchy-Born continuum model \cite{E:2007a, OrtnerTheil2012}. Let $W : \R^{d \times d} \to \R$ be a strain energy density function, the Cauchy-Born energy density $W$ is defined by
\begin{displaymath}
  W(\mF) := \det \mA^{-1} V(\mF \cdot \Rg).
\end{displaymath}

\subsection{Generic Formulation of Energy Based Atomistic/Continuum Coupling}
\label{sec:formulation:ac}
We give a generic formulation of the a/c coupling method and employ concepts and notation from various earlier works, such as \cite{Ortiz:1995a,Shenoy:1999a,Shimokawa:2004,2012-optbqce, COLZ2013}, and we adapt the formulation to the settings in this paper. 

\begin{figure}[htb]
\begin{center}
	\includegraphics[scale=0.7]{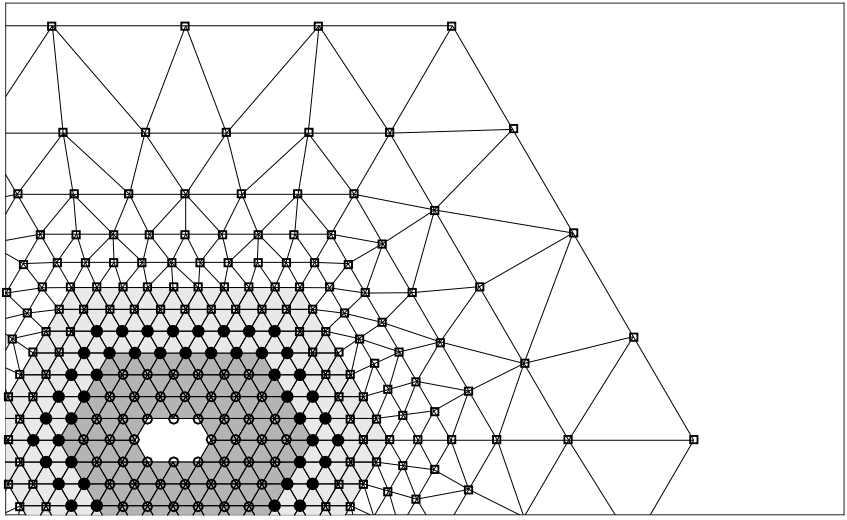}
	\caption{Illustration of computational mesh. The computational domain is $\Omega_R$, and the corresponding triangulation is $\T_{h,R}$. $\circ$ in the DimGrey region are atoms of $\L^{a}$. For the next nearest neighbour interaction, $\L^\i$ contains $\bullet$ in the LightGrey interface region. $\Box$ are continuum degrees of freedom. }
	\label{figs:plotMesh}
\end{center}
\end{figure}

The computational domain $\Omega_R = \Omega_R^\a \bigcup \Omega_R^\c \subset \R^d$ is a simply connected, polygonal and closed set, consists of the atomistic region $\Omega_R^{\a}$ and the continuum partition $\Omega_R^\c$, where $R$ is the radius of $\Omega_R$. Given the reference lattice $\Lambda$ with some local defects, we decompose the set  $\L^{\a,\i} := \L \bigcap \Omega^\a_R = \L^\a \bigcup \L^\i$ into a core atomistic set $\L^\a$ and an interfacial atomistic set $\L^\i$ such that $\L\bigcap\Ddef \subset \L^\a$, where $\Ddef$ represents the defect core. Let $\T^\a_{h,R}$ (respectively $\T^\i_{h,R}$) be the canonical triangulation induced by $\L^{\a}$ (respectively $\L^{\i}$), and $\T^{c}_{h,R}$ be a shape-regular simplicial partition of the continuum region. We denote $\Thr = \T^\c_{h,R} \bigcup \T^\i_{h,R} \bigcup \T^\a_{h,R}$ as the triangulation of the a/c coupling configuration. Please see Figure \ref{figs:plotMesh} for an illustration of the computational mesh.

The space of \textit{coarse-grained} displacements is,

\begin{align*}
  \Us_{h,R} := \b\{ u_h : \Omega_{h, R} \to \R^m \bsep ~&
  \text{ $u_h$ is continuous and p.w. affine w.r.t. $\T_{h,R}$, } \\[-1mm]
  & \text{ $u_h = 0$ on $\partial \Omega_R$ } \b\}.
\end{align*}

The subscript $R$ in the above definitions can be dropped if there is no confusion, for example, we can replace $\Thr$ by $\Th$.

Let $\vor(\ell)$ represents the voronoi cell associated with $\ell$ of the homogeneous reference lattice $\Lhom := \mA \Z^d$, for some given non-singular $\mA \in \R^{d \times d}$. We have $|\vor(\ell)| = \det \mA$, for each $\ell \in \Lhom$. For each $\ell \in \L$ denote its effective cell as $\nu_\ell$  (see \cite{PRE-ac.2dcorners}), let $\omega_\ell:=\displaystyle{\frac{|\nu_\ell|}{|\vor(\ell)|}}$ be the effective volume associated with $\ell$. For each element $T \in \Th$ we define the effective volume of $T$ by 

\begin{displaymath}
\omega_T := |T \setminus (\bigcup_{\ell \in \L^{\a}} {\textrm{vor}}(\ell))\setminus(\bigcup_{\ell \in \L^{\i}} \nu^\i_\ell)|.
\end{displaymath}

We note that $\omega_T =0 $ if $T\in \T^\a_h\setminus\T^\i_h$, $\omega_T= |T|$ if  $T\in \T^\c_h\setminus \T^\i_h$, and $0\leq \omega_T < |T|$ if $T\in \T^\i_h$, dependes on how we defining $\nu_{\ell}^{i}$, the effective cell of $\ell\in\L^{i}$.
The choices of $\nu_\ell$ and $\omega_T$ satisfy $\sum_{\ell\in \L^{\a,\i}} \nu_\ell + \sum_{T\in \Th} \omega_T = |\Omega_{h,R}|$.
	
Now we are ready to define the generic a/c coupling energy functional $\Eh$,

\begin{align}
  \label{eq:generic_ac_energy}
  \Eh(y_h) := & \sum_{ \ell \in \L^\a}  V_\ell(y_h)  + \sum_{\ell \in \L^\i} \omega_\ell V^\i_\ell(y_h)  + \sum_{T \in \Th} \omega_T  W(\D y_h|_T)   
\end{align}
where $V_\ell^\i$ is a modified interface site potential. 

The goal of a/c coupling is to find
\begin{equation}
  \label{eq:min_ac}
  y_{h,R} \in \arg\min \b\{ \Eh(y_h) \bsep y_h - y^B \in
  \Us_{h,R} \b\}.
\end{equation}

The subscript $R$ in $y_{h,R}$ and $\Us_{h,R}$ can be omitted if there is no confusion.

\def\Eh{\E^{\textrm{h}}}
The first variation of the a/c coupling variational problem \eqref{eq:min_ac} is to find $y_h-y^\mB \in \Us_{h, R}$ such that
\begin{equation}
\label{eqn:firstvariationeh}
	\<\del\Eh(y_h), v_h\> = 0, \quad\forall v_h\in \Us_{h, R}.
\end{equation}

Spurious artificial force could occur at the interface for energy based coupling even for homogeneous deformation \cite{Shenoy:1999a}, and was dubbed "ghost force". The issue of "ghost force removal" has received considerable attention in the recent years, and consistent a/c coupling methods without ghost force were developed by \cite{Shimokawa:2004,E:2006} in one dimension and \cite{PRE-ac.2dcorners,Shapeev:2010a} in two dimensions. We will introduce the consistent GRAC formulation in  \S~\ref{sec:grac} .

\def\yB{y^\mB}
\section{General GRAC formulation}
\label{sec:grac}

In this section, we describe the construction of the \textit{geometric reconstruction based consistent atomistic/continuum} (GRAC) coupling energy for multibody potentials with general interaction range and arbitrary interfaces. 

Given the homogeneous site potential $V\b( D y(\ell) \b)$, we can represent the interface potential $V_\ell^\i$ in \eqref{eq:generic_ac_energy} in terms of $V$.  For each $\ell \in
\L^\i, \rho, \vsig \in \Rg_\ell$, let $C_{\ell;\rho,\vsig}$ be
free parameters, and define
\begin{equation}
  \label{eq:defn_Phi_int}
  V_\ell^\i(y) := V \B( \b( {\textstyle \sum_{\vsig \in
      \Rg_\ell} C_{\ell;\rho,\vsig} D_\vsig y(\ell) } \b)_{\rho \in
    \Rg_\ell} \B)
\end{equation}

A convenient short-hand notation is

\begin{displaymath}
  V_\ell^\i(y) = V( C_\ell \cdot Dy(\ell) ), \quad \text{where}
  \quad \cases{
    C_\ell := (C_{\ell;\rho,\vsig})_{\rho,\vsig \in \Rg_\ell}, \quad
    \text{and} &\\
    C_\ell \cdot Dy := \b( {\textstyle \sum_{\vsig \in
      \Rg_\ell} C_{\ell;\rho,\vsig} D_\vsig y } \b)_{\rho \in
    \Rg_\ell}. & }
\end{displaymath}

We call the parameters $C_{\ell;\rho,\vsig}$ as the \textit{reconstruction parameters}. 

To construct consistent a/c coupling energy, we need to enforce the so-called \textit{patch tests} for the energy functional $\Eh$, namely, energy patch test \eqref{eq:energy_pt} and force patch test \eqref{eq:force_pt}. Those patch tests in turn prescribe conditions \eqref{eq:E_pt_grac} and \eqref{eq:F_pt_grac} for the reconstruction parameters $C$. In general, the reconstruction parameters satisfying patch tests are not unique, an $\ell^{1}$-minimization technique can be introduced to choose the "optimal" parameters \cite{COLZ2013}. Also, a stabilisation mechanism can be applied to improve the stability of the GRAC coupling scheme \cite{2013-stab.ac}.

 \subsection{Energy patch test} 
 To guarantee that $\Eh$ approximates the atomistic
energy $\E=\sum_{\ell\in\L}V_{\ell}(y)$, it is reasonable to require that the interface
potentials satisfy an \textit{energy patch test}
\begin{equation}
  \label{eq:energy_pt}
  V_\ell^\i(y^{\mF}) = V(y^{\mF}) \qquad \quad \forall\ell \in \L^\i, \quad \mF \in \R^{m
    \times d}.
\end{equation}
namely, the interface potential coincides with the atomistic potential for uniform deformations.
 
For the GRAC coupling scheme, a sufficient and necessary condition for the energy patch test is that $\mF \cdot \Rg(\ell) = C_\ell \cdot (\mF \cdot \Rg)$ for all $\mF \in \R^{m \times d}$ and $\ell \in \L^\i$. This is equivalent to
\begin{equation}
  \label{eq:E_pt_grac}
  \rho = \sum_{\vsig \in \Rg(\ell)} C_{\ell; \rho,\vsig} \vsig \qquad
  \forall \ell \in \L^\i, \quad \rho \in \Rg(\ell).
\end{equation}

\subsection{Force patch test}
\label{sec:force_pt_lineqn}
We call the following condition the \textit{force patch test}, namely, for $\L = \L^{\textrm{hom}}$ and $\Phi_{\ell} = \Phi$,
 
\begin{equation}
  \label{eq:force_pt}
  \< \del \Eh(y^{\mF}), v_h \> = 0 \qquad \forall v_h \in \UsT,
  \quad \mF \in \R^{m \times d}.
\end{equation}
where $\mF$ is some uniform deformation gradient. This is saying that there is no artificial "ghost force" for uniform deformations. 

From the general GRAC formulation \eqref{eq:generic_ac_energy}, we can decompose the first variation of the a/c coupling energy into three parts,
\begin{equation}
  \<\del\Eh(y^{\mF}), v_h\> = \<\del\Ea(y^{\mF}), v_h\> 
  + \<\del\Ei(y^{\mF}), v_h\> + \<\del\Ec(y^{\mF}), v_h\>.\quad \forall v_h \in \UsT
  \label{eq:threeparts}
\end{equation}

To simplify the notation, we drop the $y^{\mF}$ dependence from the following expressions in this section, for example, we write $\Ea$ instead of $\Ea(y^{\mF})$,
$\D_\rho V$ instead of $\D_\rho V(Dy^{\mF})$, and so forth. Here, $\D_\rho
V$ denotes the partial derivative of $V$ with respect to the $D_\rho y$
component. Since $\D_\rho V = -\D_{-\rho} V$, we only consider half of the directions in the
interaction range: fix $\Rgp \subset \Rg$ such that $\Rgp \cup
(-\Rgp) = \Rg$ and $\Rgp \cap (-\Rgp) = \emptyset$.

As proposed in \cite{COLZ2013}, a necessary and sufficient condition on the reconstruction parameters $C_\ell$ to satisfy the force patch test
  \eqref{eq:force_pt}  for all $V \in C^\infty((\R^d)^\Rg)$ is
  \begin{equation}
    \label{eq:F_pt_grac}
    c^\a_\rho(\ell) + c^\i_\rho(\ell) + c^\c_\rho(\ell) = 0,
  \end{equation}
  for $\ell\in\Li+\Rg$, and $\rho\in\Rgp$. The coefficients $c^\a_\rho(\ell)$, {$c^\i_\rho(\ell)$} and
$c^\c_\rho(\ell)$ are geometric parameters with respect to the underlying lattice and the interface geometry, formulated by collecting all the coefficients for the terms $\D_\rho V \cdot$ of the first variations of a/c coupling energy as in the following equations \eqref{eq:firstvar}. The interface coefficients $c^\i_\rho(\ell)$ depend linearly on the unknown reconstruction paramters $C_{\ell;\rho,\vsig}$. Then, by \eqref{eq:threeparts}, we have
   
\begin{align}
	\begin{split}
  \<\del\Ea, v_h\>  & = \sum_{\ell\in\La +\Rg}\sum_{\rho\in\Rgp}
        c^\a_\rho(\ell) \b[\D_\rho V \cdot v_h(\ell)\b] \\
   &= \sum_{\substack{\rho \in \Rgp \\ \ell \in
      \La-\rho}} {\b[\D_\rho V \cdot v_h(\ell) \b]}- \sum_{\substack{\rho \in
      \Rgp \\ \ell\in\La+\rho}}\b[ \D_\rho V \cdot v_h(\ell)\b], \\
  \<\del\Ei, v_h\>  & = \sum_{\ell\in\L^i
          +\Rg}\sum_{\rho\in\Rgp}c^\i_\rho(\ell) \b[ \D_\rho V \cdot v_h(\ell) \b] \\
   &= \sum_{\substack{\vsig \in \Rg \\ \ell \in
      \Li+\vsig}} \omega_{\ell-\vsig}^\i \sum_{\rho\in\Rgp}
  (C_{\ell-\vsig;\rho,\vsig}-C_{\ell-\vsig;-\rho,\vsig}) \b[\D_\rho V \cdot v_h(\ell)\b]\\
  & \qquad
  -\sum_{\ell\in\Li} \omega_\ell^\i \sum_{\rho\in\Rgp}
  \sum_{\vsig\in\Rg} (C_{\ell;\rho,\vsig}-C_{\ell;-\rho,\vsig}) \b[\D_\rho
  V \cdot v_h(\ell) \b], \quad \text{and} \\
  \<\del\Ec, v_h\>  &= \sum_T\sum_{\rho\in\Rgp}\sum_{i=1}^3
  2\frac{\omega_T}{\det\mA}\D_T\phi_i^T\cdot \rho \b[ {\D_\rho V} \cdot v^T_{h,i} \b], \\
  &= \sum_T\sum_{\rho\in\Rgp}\sum_{i=1}^3
  2\frac{\omega_T}{\det\mA}\D_T\phi_i^T\cdot \rho \b[ {\D_\rho V} \cdot
  v^T_{h,i} \b],
  \end{split}
  \label{eq:firstvar}
\end{align}
where the nodes $\ell_{i}^{T}$ are the three corners of the triangle $T$, $v_{h,i}^T = v(\ell_i^T)$ and $\phi_i^T$ are the three nodal linear bases corresponding to $v_{h,i}^T$, $i = 1,2,3$, 

The force patch test is automatically satisfied for the atomistic model and the Cauchy-Born continuum model.Therefore, we only need to consider the force consistency for those sites with the extended interface region $\L^\i+\Rg:=\{\ell\in\L|\exists \ell'\in\L^\i, \exists\rho\in\Rg,\text{ such that } \ell = \ell'+\rho\}$. Therefore \eqref{eq:F_pt_grac} for $\ell\in\L^\i+\Rg$ together with \eqref{eq:E_pt_grac} for $\ell\in\L^\i$ form a linear system for the unknown parameters $C_{\ell;\rho,\vsig}$. 

\begin{remark}
In \cite{PRE-ac.2dcorners,COLZ2013}, we choose $ \omega_\ell^i = {\vor}(\ell), \forall \ell \in \L^i$. To reduce the degrees of freedom of the reconstruction parameters and also the number of constraint equations, we can use a variation of the local reflection method in \cite{COLZ2013}. Namely,  we choose $\omega_\ell^i  = {\vor}(\ell) \cap \Omega^\a$, and enforce that $C_{\ell; \rho,\vsig} = 0$ for  $\ell \in \Li$ and $\ell+\vsig\in \Omega^\c$, then we only need to impose the force balance equation for $(\Li + \Rg) \cap \L^{\a}$. It has been shown that in \cite{COLZ2013}, the linear system for the reconstruction parameters is underdetermined, and up to numerical accuracy, the solution exists and therefore it is not unique.

\label{rm:rank_def}
\end{remark}

\begin{figure}
  \begin{center}
    \includegraphics[width=6cm]{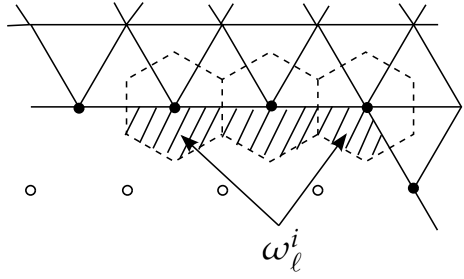}
    \caption{Effective Voronoi cells for the
      interface nodes (filled circles) are the shaded area in the
      above figure. $\omega_\ell^\i<1$ for the outmost interface atoms which are
      adjacent to the continuum region.}\label{fig:effvol}
  \end{center}
\end{figure}


%

\subsection{Consistency and Optimisation of Reconstruction Parameters}
\label{sec:optim_coeffs}
Due to the non-uniqueness of the reconstruction parameters $C_{\ell,\rho,\vsig}$, we need to choose the "optimal" parameters. A naive idea is to use a least-squares approach to minimize $\|C\|_{\ell^2}$,
\begin{equation}
  \label{eq:least_squares}
  \text{minimize } \sum_{\ell \in \L^\i} \sum_{\rho,\vsig \in
    \Rg(\ell)} |C_{\ell;\rho,\vsig}|^2 \quad \text{ subject to \eqref{eq:E_pt_grac} and \eqref{eq:F_pt_grac}.}
\end{equation}
However, the resulting parameters do not give a convergent method.

It is shown in \cite[Thm. 6.1]{Or:2011a} , under the assumptions
that $d = 2$ and that the atomistic region $\Omega^\a$ is connected, any a/c coupling
scheme of the type \eqref{eq:generic_ac_energy} satisfying the energy
and force patch tests \eqref{eq:energy_pt},  \eqref{eq:force_pt} is first-order consistent: if $y = \yB$ in
$\L \setminus \Omega_R$ and if $\tilde{y}$ is an $H^2_{\textrm
  {loc}}$-conforming interpolant of $y$, then
\begin{equation}
  \label{eq:consistency}
  \b\< \del \E(y) - \del \Eh(I_h y), v_h \b\> \leq C_1
  \| h \D^2 \tilde{y} \|_{L^2(\tilde{\Omega}^\c)} \|\nabla v_h\|_{L^2},
\end{equation}
where $C_1$ is independent of $y$. The dependence of $C_1$
on the reconstruction parameters $C_\ell$ is analysed in \cite{COLZ2013},
\begin{align}
  \label{eq:C1_dependence}
  & C_1 \leq C_1' \, (1 + \textrm{width}(\L^\i)) \, \sum_{\rho,\vsig \in
    \Rg} |\rho|\,|\vsig|\, M_{\rho,\vsig} + C_1''\\
  \notag
  \notag & \text{where} \quad M_{\rho,\vsig} = \max_{\ell \in \L^\i}
  \sum_{\tau,\tau' \in \Rg(\ell)} |V_{\tau,\tau'}(C_\ell \cdot Dy(\ell))| \,
  |C_{\ell;\tau,\rho}|\,|C_{\ell;\tau',\vsig}|.
\end{align}
$C_1'$ is a generic constant and $C_1''$ does not depend on the
reconstruction parameters.

%

Intuitively one may think of 
\begin{displaymath}
  M(\ell) := \sum_{\rho,\vsig} |\rho|\,|\vsig| \sum_{\tau,\tau'}
  \b|V_{\tau,\tau'}(C_\ell \cdot Dy(\ell))\b|  \,
  |C_{\ell;\tau,\rho}|\,|C_{\ell;\tau',\vsig}|
\end{displaymath}
to be a realistic ($\ell$-dependent) pre-factor. With generic structural assumption $|V_{\tau,\tau'}(C_\ell \cdot Dy(\ell))| \lesssim
\omega(|\tau|)\,\omega(|\tau'|)$, where $\omega$ has some decay that
is determined by the specific interaction potential, see the discussion for an EAM type potential in App.B.2 \cite{OrtnerTheil2012}. We obtain that

\begin{align*}
  M(\ell) &\lesssim \sum_{\rho,\vsig} |\rho|\,|\vsig| \sum_{\tau,\tau'}
  \omega(|\tau|)\,\omega(|\tau'|) \,
  |C_{\ell;\tau,\rho}|\,|C_{\ell;\tau',\vsig}| \\
  &= \B(\sum_{\rho,\tau} |\rho| \omega(|\tau|)
  |C_{\ell;\tau,\rho}|\B)\,
  \B(\sum_{\vsig,\tau'} |\vsig| \omega(|\tau'|)
  |C_{\ell;\tau',\vsig}|\B) \\
  &= \B(\sum_{\rho,\tau} |\rho| \omega(|\tau|)
  |C_{\ell;\tau,\rho}|\B)^2.
\end{align*}

This indicates that, instead of $\|C\|_{\ell^2}$, we should minimize
$\max_{\ell \in \L^\i} \sum_{\rho,\tau} |\rho| \omega(|\tau|)
|C_{\ell;\tau,\rho}|$.  Since we do not in general know the generic
weights $\omega$, we simply drop them, and instead minimise
$\sum_{\rho,\tau} |C_{\ell;\tau,\rho}|$. Further, taking the maximum
of $\ell \in \L^\i$ leads to a difficult and computationally expensive
multi-objective optimisation problem. Instead, we propose to minimise
the $\ell^1$-norm of all the coefficients:
\begin{equation}
  \label{eq:ell1-min}
  \text{ minimise } \sum_{\ell \in \L^\i} \sum_{\rho,\vsig \in
    \Rg(\ell)} |C_{\ell;\rho,\vsig}| \quad \text{ subject to
    \eqref{eq:E_pt_grac} and \eqref{eq:F_pt_grac}.} 
\end{equation}

We remark that with the simplifications, the reconstruction parameters only depend on the interaction range, not the particular form of interaction potentials, also, $\ell^1$-minimisation tends to generate ``sparse'' reconstruction
parameters which may present some gain in computational cost in the
energy and force assembly routines for $E^\ac$.

\subsection{Stability and stabilisation}
\label{sec:stab}

The stability of a/c coupling method is of great importance for numerical analysis and the issue of critical stability is essential for practitioners.
The issue of stability and (in)stability of a/c coupling scheme is discussed in detail in \cite{2013-stab.ac}, in particular, a stabilization mechanism
is proposed to reduce the stability gap between the a/c coupling methods and the corresponding atomistic model (ground truth). We sketch those results in this subsection. 

We expect a/c coupling method has the following stability estimate of the form,
\begin{equation}
  \label{eq:stability}
  \< \ddel \Eh(I_h y) v_h, v_h \> 
  \geq c_0 \| \nabla v_h \|_{L^2}^2
\end{equation}
together with consistency and other technical assumptions, we can prove the convergence of the a/c coupling method with inverse function theorem with the form $\|u - u_h\| \lesssim \text{consistency error} / c_0$. 

Also, we have the observation that for any uniform deformation, the stability constant of a/c coupling method is less or equal to the stability constant of the corresponding atomistic model. This motivate the introduction of the following stabilization of the interface potential, 

\begin{equation}
  \label{eq:qnl_stab}
  V_\ell^\i(y_h) := 
  V \b( C_\ell \cdot Dy_h(\ell) \b) + \kappa |D_{\textrm{nn}}^2 y_h(\ell)|^2,
\end{equation}
where $\kappa \geq 0$ is a stabilisation parameter, and $|D_{\textrm{nn}}^2
y_h(\ell)|^2$ is defined as follows: we choose $m \geq d$ linearly
independent ``nearest-neighbour'' directions $a_1, \dots, a_m$ in the
lattice, and denote
\begin{displaymath}
  \b|D_{\textrm{nn}}^2 y_h(\ell)\b|^2 := \sum_{j = 1}^m \b| y_h(\ell+a_j) -
  2 y_h(\ell) + y_h(\ell-a_j) \b|^2.
\end{displaymath}

If $\kappa$ is $O(1)$ constant, the stabilization term does not affect consistency, the reconstruction parameters $C_\ell$ are still determined by \eqref{eq:ell1-min}. However, the stabilisation will affect the computation of the stress tensor in the interface region, which we will introduce in the following section.
In \cite{2013-stab.ac}, we have theoretically proved for some symmetric configuration and numerical justified for prototypical examples that such stabilization can reduce the stability gap and suppress spurious critical mode for the a/c coupling methods. 


\section{Stress formulation}
\label{sec:stress}
\def\Ta{\T_\a}
\def\Th{\T_{\textrm{h}}}
\def\sh{\sigma^{\textrm{h}}}
\def\si{\sigma^{\textrm{i}}}


The stress formulation plays a significant role in the analysis of a/c coupling methods. Given deformation $y-y^B\in \Use$ and an energy functional $\E(y)$, the first Piola-Kirchhoff stress $\sigma(y)$ satisfies 
\begin{displaymath}
  \<\delta\E(y), v\> = \int_{\R^d}\s(y):\D v \dx,
\end{displaymath}
for any $v \in \Use$. 

From the first variation of the a/c energy functional \eqref{eq:threeparts} and \eqref{eq:firstvar}, we expect to derive its stress formulation, namely, for any $v_h\in\Ush$,
\begin{align*}
 \<\delta\Eh(y_h), v_h\> = &  \<\delta\Ea, v_h\> + \<\delta\Ei, v_h\> + \<\delta\Ec, v_h\> \\
 = & \int_{\Omega^\a}\sa:\D v_h \dx + \int_{\Omega^\i}\si:\D v_h \dx + \int_{\Omega^\c}\sc:\D_{T} v_h \\
 = & \int_{\Omega} \sh : \D v_h \dx 
\end{align*}

 In this section, we first review the stress tensor formula derived by \cite{OrtnerTheil2012} in \S~\ref{sec:stress:apriori} which is essential for the a priori analysis, however, this formula is not suitable for the a posteriori estimate. In \S~\ref{sec:stress:apost} we introduce a novel computable stress tensor expression, which is convenient for the purpose of a posterior error estimation and adaptive algorithm. We discuss the assembly of stress tensor for our model problem in \S~\ref{sec:stress:assm}.




\subsection{Stress Tensor Formulation in \cite{OrtnerTheil2012}}
\label{sec:stress:apriori}

We first discuss the formula for $\sa$. For simplicity, consider full atomistic energy $\E$ in \eqref{eqn:Ea}. The "canonical weak form" of $\delta\E$ is 
\begin{equation}
	\<\delta \E(u), v\> = \sum_{\ell\in \L}\sum_{\rho\in\Rg} V_{\ell, \rho}(u)\cdot D_\rho v(\ell), \quad \text{for } v\in \Use
	\label{eq:weakform}
\end{equation}

Now we define a modified version of the canonical weak form. Let 
\begin{equation}
  v^{\ast} := \bar{\zeta} \ast \bar{v}.
  \label{eq:vast}
\end{equation}

The finite differences $D_{\rho}v^{\ast}(\zeta)$ can be expressed as in \cite{Shapeev2012}, 
\begin{align}
\nonumber D_{\rho}v^{\ast}(\ell) = & \int_{s=0}^{1}\D_{\rho}v^{\ast}(\ell+s\rho) \ds = \int_{\R^d}\int_{s=0}^{1}\bar{\zeta}(\ell+s\rho-x)\D_{\rho}\bar{v}(x)\ds\dx \\
=& \int_{\R^d} \chi_{\ell,\rho}(x)\D_{\rho}\bar{v}(x)\dx
 \label{eq:dvast}
\end{align}
where $\chi_{\ell,\rho}$ is a generic weighting function defined as below, can be understood as a mollified version of the line measure.

\begin{equation}
\label{eq:chi}
\chi_{\ell,\rho}(x) := \int_{0}^{1}\bar{\zeta}(\ell+t\rho-x)\dt
\end{equation}

Now, we replace the test function in \eqref{eq:weakform} from $v$ to $v^{\ast}$, 
\begin{align*}
 \<\delta\E(y), v^{\ast}\> = & \sum_{\ell\in\L^{a}}\sum_{\rho\in\Rg}\partial_{\rho}V_{\ell}\cdot\int_{\R^d} \chi_{\ell,\rho}(x)\D_{\rho}\bar{v}(x)\dx \\
  =& \int_{\R^d}\left\{ \sum_{\ell\in\L^{a}}\sum_{\rho\in\Rg}\left[\partial_{\rho}V_{\ell}\otimes\rho\right]\chi_{\ell,\rho}(x) \right\}:\D\bar{v}\dx
\end{align*}

Thus, we have shown that for $y-y^{\mB}\in\Use$ and $v\in\Us$ with compact support , 
\begin{displaymath}
  \<\delta\E(y), v^{\ast}\> = \int_{\R^d}\sa(y;x):\D\bar{v} \dx,
\end{displaymath}
with
\begin{equation}
\label{eq:sa:priori}
\sa(y;x):=\sum_{\ell\in\L^\a}\sum_{\rho\in\Rg}\left[\partial_{\rho}V_{\ell}\otimes\rho\right]\chi_{\ell,\rho}(x),
\end{equation}

Through an analogy to the analysis above, $\<\delta\Ei(y), v^{\ast}\>$ could be written as,
\begin{displaymath}
  \<\delta\Ei(y), v^{\ast}\> = \int_{\R^d}\si(y;x):\D\bar{v} \dx,
\end{displaymath}
where
\begin{equation}
\label{eq:si:priori}
\si(y;x):=\sum_{\ell\in\L^\i}\sum_{\rho\in\Rg}\left[\omega^\i_\ell\partial_{\rho}V^{\i}_{\ell}\otimes\rho\right]\chi_{\ell,\rho}(x),
\end{equation}

Finally, the first Piola-Kirchhoff stress of the Cauchy-Born model gives,
\begin{equation}
\label{eq:sc:priori}
\sc(y;x)=\partial W(\D y(x))
\end{equation}

By proper regularity assumptions on $V$ and $y$, it can be shown that $\sc(y;x)$ is second order consistent to $\sa(y;x)$.

%
%

The stress tensor expression introduced in this section is important for the a priori analysis of Caucy-Born continuum model in \cite{OrtnerTheil2012} and blended a/c coupling  method \cite{LiOrShVK:2014}. However, since $\sa(y;x)$ and $\si(y;x)$ depend on $x$ through $\chi_{\ell,\rho}(x)$, it is relatively difficult to calculate the value of them, which is crucial for the a posteriori error estimates and adaptive algorithm. In the next section, we will introduce a computable expression of stress tensor.

\subsection{A Computable Stress Tensor Formulation}
\label{sec:stress:apost}
\def\meas{\textrm{meas}}
\def\Tlrho{\T_\ell^\rho}
\def\wlrho{\omega_{\ell}^\rho}

For the nearest neighbour interactions \cite{APEst:2017}, we use the canonical expressions of $\delta \E$ ($\delta\Eh$) to define $\sa$ ($\sh$). In that case, $\sa$ ($\sh$) is piecewise constant over the triangulation $\Ta$ ($\Th$). In this section, we will extend this formulation to general finite range interactions. 

\def\length{\textrm{length}}

Consider the canonical weak form of $\delta\E$ in \eqref{eq:weakform}, instead of replacing $v$ with $v^{\ast}$, we distribute $D_\rho v$ to relevant triangles in order to transfer the sum with respect to atoms to the sum over elements in the micro-triangulation $\T_\a$ induced by the reference lattice $\L$. We express $D_\rho v$ as,
\begin{equation}
  D_{\rho}v(\ell) = \sum_{T\in\Tlrho}\wlrho (T)\D_T v\cdot\rho
  \label{eq:stress:drho}
\end{equation}
where $\T_\ell^\rho := \{T\in \Ta|\length(T\cap(\ell, \ell+\rho))>0\}$ is the set which contains the elements that form the compact support of bond $(\ell, \ell+\rho)$, $\wlrho (T)$ is an appropriate weight function. The implementation details for two dimensional triangular lattice will be discussed in detail in \S~\ref{sec:stress:assm}. 

For any $y-y^\mB\in\Use, v\in\Use$,
\begin{equation}
\begin{split}
  \<\del\E(y), v\>&=\sum_{\ell \in \L}\sum_{\rho \in \Rg_\ell}\partial_{\rho}V_{\ell}\cdot D_{\rho}v  \\
                           &=\sum_{\ell \in \L}\sum_{\rho \in \Rg_\ell}\partial_{\rho}V_{\ell}\cdot\left(\omega_{\rho}\left(\sum_{b=(\ell, \ell+\rho)\cap  T \neq \emptyset}\Dc{T} v\right)\cdot\rho\right) \\
                           &=\sum_{\ell \in \L}\sum_{\rho \in \Rg_\ell}\frac{2|T|}{\det\mA}\partial_{\rho}V_{\ell}\otimes\rho:\left(\omega_{\rho}\sum_{b=(\ell, \ell+\rho)\cap  T \neq \emptyset}\Dc{T} v\right)  \\
                           &=|T|\sum_{T\in \Ta}\sum_{\rho \in \Rg_\ell}\frac{1}{\det \mA}\sum_{b=(\ell, \ell+\rho)\cap  T \neq \emptyset}2\omega_{\rho}\partial_{\rho}V_{\ell}\otimes\rho:\left(\Dc{T} v\right), \forall v\in\Us  
\end{split}
\label{eq:atomvar}
\end{equation}
$\wlrho(T)$ represents the contribution of the value of $\partial_{\rho}V(\Dc{T}y_h)\otimes\rho$ from element $T$, which depends on the specific type of $\rho\in\Rg$. 


Recall from \eqref{eq:generic_ac_energy} and \eqref{eq:defn_Phi_int}, the GRAC a/c coupling energy functional is of the form

\begin{displaymath}
  \Eh(y_h) :=  \sum_{ \ell \in \L^\a}  V_\ell(y_h)  + \sum_{\ell \in \L^\i} \omega_\ell V^\i_\ell(y_h)  + \sum_{T \in \Th} \omega_T  W(\D y_h|_T)   
\end{displaymath} 
with 
\begin{displaymath}
  V_\ell^\i(y) := V_{\ell} \B( \b( { \sum_{\vsig \in
      \Rg_\ell} C_{\ell;\rho,\vsig} D_\vsig y(\ell) } \b)_{\rho \in
    \Rg_\ell} \B)
\end{displaymath}

With the decomposition of the first variation of GRAC in \eqref{eq:threeparts}, the first variation of the interface energy $\Ei$ has the following form, 
\begin{equation}
	\begin{split}
  \<\del\sum_{\ell \in \L^{\i}} \omega_{\ell} V^{\i}_{\ell}(y_h), v_{h}\>  &= \<\del\sum_{\ell \in \L^{\i}} \omega_{\ell} V_{\ell} \B( \b( { \sum_{\vsig \in\Rg_\ell} C_{\ell;\rho,\vsig} D_{\vsig} y_h } \b)_{\rho \in\Rg_\ell} \B), v_{h}\>  \\
        &= \sum_{\ell \in \L^\i} \omega_\ell \sum_{\rho \in\Rg_\ell}  \partial_{\rho}V_{\ell} \cdot \left(\sum_{\vsig \in\Rg_\ell} C_{\ell;\rho,\vsig} D_\vsig v_h(\ell)\right) \\
        &= \sum_{\ell \in \L^\i} \omega_\ell \sum_{\rho \in\Rg_\ell} \sum_{\vsig \in\Rg_\ell} \partial_{\rho}V_{\ell} C_{\ell;\rho,\vsig}  \cdot \left(\omega_{\vsig}\left(\sum_{b=(\ell, \ell+\vsig)\cap  T \neq \emptyset}\Dc{T} v_{h}\right)\cdot\vsig \right)   \\
        &= \sum_{\ell \in \L^\i} \frac{2\omega_\ell |T|}{\det\mA} \sum_{\vsig \in\Rg_\ell} \left(\sum_{\rho \in\Rg_\ell}C_{\ell;\rho,\vsig} \partial_{\rho}V_{\ell} \otimes \vsig\right)  : \left(\omega_{\vsig}\sum_{b=(\ell, \ell+\vsig)\cap  T \neq \emptyset}\Dc{T} v_{h}\right)  \\
        &= \sum_{T\in\T_{\h}^{\i}} \frac{|T|}{\det\mA} \sum_{\vsig \in\Rg_\ell} \sum_{b=(\ell, \ell+\vsig)\cap  T \neq \emptyset} 2\omega_{\ell}\left(\sum_{\rho \in\Rg_\ell}C_{\ell;\rho,\vsig} \omega_{\vsig}\partial_{\rho}V_{\ell} \otimes \vsig\right)  : \left(\Dc{T} v_{h}\right) 
        \end{split}
        \label{eq:intfvar}
\end{equation}


Consider the continuum model, for any $y_h - y^\mB\in \Ush,v_h\in\Ush$, we have
\begin{align}
  \<\del\Ec(y_h),v_h\> &=\sum_{T\in\Th}|T|\pp W(\Dc{T} y_h)  \label{eq:contvar} \\
  				&= \sum_{T\in\Th}|T|\sum_{\rho\in\Rg} \partial_{\rho} V(\nabla_T y_h) \otimes \rho:\Dc{T} v_{h}, 
\end{align} 


Adding the stabilization term as in \eqref{eq:qnl_stab}, we are ready to obtain the first variation of the (stabilized) GRAC coupling method,

\def\Rnn{\R_{\textrm{nn}}}
\begin{align}
    \<\del\Eh(y_h),v_h\>&=\sum_{T\in\T_{\h}^{\a}\cup\T_{\h}^{\i}} \frac{|T|}{\det \mA} \sum_{\rho \in\Rg_\ell} \sum_{b=(\ell, \ell+\rho)\cap  T \neq \emptyset}2\omega_{\ell} \wlrho(T)\partial_{\rho}V_{\ell}^{h} \otimes \rho  : \left(\Dc{T} v_{h}\right)  
  \nonumber \\
  & + \sum_{T\in\Th} \frac{\omega_T}{\det \mA}\sum_{\rho\in\Rg} \partial_{\rho} V(\nabla_T y_h) \otimes \rho:\Dc{T} v_{h}, \label{eq:acvar} \\
  & + C_{\textrm{stab}}\sum_{T'\in\Th^{i}}\sum_{\zeta\in\Rnn}\sum_{b'=(\ell, \ell+\zeta)\cap  T' \neq \emptyset}2\omega_{\ell}^{\zeta}(T) \textrm{d}_{\textrm{nn}}^{3}(\ell, \zeta)\otimes\zeta : \Dc{T'} \nonumber
\end{align}

\begin{displaymath}
  \text{where}
  \quad \partial_{\rho}V_{\ell}^{h} := \cases{
    \partial_{\rho}V_{\ell}, \quad \ell\in\L^{a}
     &\\
    \sum_{\vsig \in\Rg_\ell}C_{\ell;\vsig,\rho} \partial_{\vsig}V_{\ell}, \quad \ell\in\L^{i} & }
\end{displaymath}
and
\begin{displaymath}
 \textrm{d}_{\textrm{nn}}^{3}(\ell, \zeta) = -y_h(\ell+2\zeta)+3y_h(\ell+\zeta)-3y_h(\ell)+y_h(\ell-\zeta)
\end{displaymath}

$C_{\textrm{stab}}$ is a constant equal to 1 if stabilisation is applied, and 0 otherwise.

We can define the \textit{atomistic stress tensor} $\sa$ (with respect to micro-triangulation $\Ta$), the \textit{continuum
stress tensor} $\sc$, and the \textit{a/c stress tensor} $\sh$ by the following first variations of different models,
\begin{align}
  \<\del\Ea(y),v\>&=\sum_{T\in\Ta}|T|\sa(y;T):\Dc{T} v, \forall v\in\Us, 
  \label{eq:atomstress}\\
  \<\del\Ec(y_h),v_h\>&=\sum_{T\in\Th}|T|\sc(y_h;T):\Dc{T} v_h, \forall v_h\in\Ush,
  \label{eq:contstress}\\
    \<\del\Eh(y_h),v_h\>&=\sum_{T\in\Th}|T|\sh(y_h;T):\Dc{T} v_h, \forall v_h\in\Ush.
  \label{eq:acstress}
\end{align}

From the above discussions, we have the following "canonical" choices for $\sa$, $\sc$ and $\sh$. 
\begin{align}
    \label{eq:defn_Sa}
    \sa(y; T) :=~&  \frac{1}{\det \mA}\sum_{\rho \in \Rg_\ell}\sum_{b=(\ell, \ell+\rho)\cap  T \neq \emptyset}2\wlrho(T)\partial_{\rho}V_{\ell}\otimes\rho, \\
    \label{eq:defn_Sc}
    \sc(y_h; T) :=~& \frac{1}{\det \mA}\sum_{\rho\in\Rg} \partial_{\rho} V(\nabla_T y_h) \otimes \rho, \\
    \label{eq:defn_Sh}
    \sh(y_h; T) :=~&\frac{1}{\det \mA}\left(\sum_{\rho \in\Rg_\ell} \sum_{b=(\ell, \ell+\rho)\cap  T \neq \emptyset} 2\omega_{\ell}\wlrho(T) \partial_{\rho}V_{\ell}^{h} \otimes \rho + \frac{\omega_T}{|T|} \sc(y_h; T) \right.  \\
    & \left. +C_{\textrm{stab}}\sum_{\zeta\in\Rnn}\sum_{\substack{b'=(\ell, \ell+\zeta)\cap  T' \neq \emptyset \\ T'\in\T^\i}} 2\omega_{\ell}^{\zeta}(T) \textrm{d}_{\textrm{nn}}^{3}(\ell,\zeta)\otimes\zeta \right). \nonumber
\end{align}

However, the stress tensors defined through \eqref{eq:atomstress}-\eqref{eq:acstress} are not unique due to the following results.

\begin{definition}
We call piecewise constant tensor field $\sigma\in \PO(\T)^{2\times 2}$ \textit{divergence free} if 
\begin{displaymath}
	\sum_{T\in\T}|T|\sigma(T):\Dc{T} v\equiv 0, \forall v\in ({\textrm{P}}_1(\T))^2.
\end{displaymath}
\end{definition}


\begin{corollary}
By definitions \eqref{eq:acstress}, it is easy to know that the force patch test condition \eqref{eq:force_pt} is equivalent to that $\sh(y_\mF)$ is divergence free for any constant deformation gradient $\mF$.
\end{corollary}


\def\Fh {\mathcal{F}_h}

 The discrete divergence free tensor fields over the triangulation $\T$ can be characterized by the non-conforming Crouzeix-Raviart finite elements \cite{PRE-ac.2dcorners, Or:2011a}. The Crouzeix-Raviart finite element space over $\T$ is defined as
 
\begin{align*}
N_{1}(\T)=\{c:\bigcup_{T\in\T}\textrm{int}(T)\to\R \quad \big{|}& \quad c \textrm{ is piecewise affine w.r.t. }\T, \textrm{and} \\ 
 & \textrm{continuous in edge midpoints } q_{f}, \forall f\in\mathcal{F}\}
\end{align*}

The following lemma in \cite{PRE-ac.2dcorners} characterizes the discrete divergence-free tensor field.

\begin{lemma}
A tensor field $\sigma\in\mP_{0}(\T)^{2\times 2}$ is divergence free if and only if there exists a constant $\sigma_{0}\in\R^{2\times 2}$ and a function $c\in N_{1}(\T)^{2}$ such that 
\begin{displaymath}
\sigma = \sigma_{0} + \nabla c\mJ, \qquad \textrm{where}\quad\mJ = \mymat{0 & -1 \\ 1 & 0}\in{\textsf{SO}}(2).
\end{displaymath}
\label{lem:divfree}
\end{lemma}

The following corollary provides a representation of the stress tensors defined in \eqref{eq:atomstress}-\eqref{eq:acstress}.

\begin{corollary}
\label{cor:divfree}
The stress tensors in the definitions \eqref{eq:atomstress}-\eqref{eq:acstress} are not unique. 
Given any stress tensor $\sigma \in \mP_{0}(\T)^{2\times 2}$ satisfies one of the definitions \eqref{eq:atomstress}-\eqref{eq:acstress} , where $\T$ is the corresponding triangulation. Define the admissible set as $\textrm{Adm}(\sigma):=\{\sigma  + \nabla c \mJ, c\in N_1(\T)^2\}$, then any $\sigma'\in \textrm{Adm}(\sigma)$  satisfies the definition of stress tensor.
\end{corollary}

\def\Pdef{\mathscr{P}^{\textrm{def}}}
\def\Usc{\Us^c}
\def\Use{\Us^{1,2}}
\def\Usrh{\Us_R^h}
\def\Ush{\Us^h}
\def\Eh{\E^{\textrm{h}}}
\def\Tr{\textrm{Tr}}
\def\sjump{\llbracket\sh\rrbracket_f}
\def\dx {\textrm{dx}}
\def\Fh {\mathcal{F}_h}
\def\Fhi{\Fh\bigcap {\textrm{int}}(\Omega_R)}

\section{Adaptive Algorithms and Numerical Experiments}
\label{sec:numerics}
\def\Ta{\T_\a}
\def\Th{\T_{\textrm{h}}}
\def\sh{\sigma^{\textrm{h}}}
\def\sjump{\llbracket\sh\rrbracket}

\newcommand{\fig}[1]{Figure \ref{#1}}
\newcommand{\tab}[1]{Table \ref{#1}}

We have derived the following a posteriori estimate in \cite[Theorem 3.1]{APEst:2017}: let $y_h$ be the a/c solution, for any $v\in\Use$, the residual $\mR[v] = \< \del\E(\Ia y_h), v \>$ is bounded by the sum of the following estimators,
\begin{equation}
	 \< \del\E(\Ia y_h), v \> \leq  \big(\eta_T(y_h) + \eta_M(y_h) + \eta_C(y_h) \big)\|\nabla v\|_{L^2}
	\label{eq:residual}
\end{equation}
where $\eta_T$ is the truncation error estimator (the $L^2$ norm of the atomistc stress tensor close to the outer boundary), $\eta_M$ is the modelling error estimator (the difference of a/c stress tensor and atomistic stress tensor), and $\eta_C$ is the coarsening error (jump of a/c stress tensor across interior edges). They are given by 
\begin{equation}
\label{eqn:etat}
	\eta_T(y_h): = C_1\|\sa(\Ia y_h)-\sigma^\mB \|_{L^2(\Omega_R\setminus B_{R/2})}.
\end{equation}
where $\displaystyle\sigma^\mB = \frac{1}{\det \mA}\sum_{\rho \in \Rg_\ell}\partial_{\rho}V(\mB a)\otimes\rho$.


\begin{equation}
\label{eqn:etam}
	\eta_M(y_h) := C_2\big\{\sum_{T\in\Ta}|T|\big[\sa(\Ia y_h, T) - \sum_{T'\in \Th, T'\bigcap T\neq \emptyset}\frac{|T'\bigcap T|}{|T|}\sh(y_h, T')\big]^2\big\}^{\frac12}.
\end{equation}

\begin{equation}
\label{eqn:etac}
\eta_C(u_h) := C_3(\sum_{f\in\Fh} (h_f\sjump)^2)^\frac12
\end{equation}
where $C_1$, $C_2$, and $C_3$ are independent of $R$ and $u_h$, actually $C_3 = \sqrt{3}CC'_{\Th}$ depends only on the shape regularity of $\Th$.

In this section, we will propose an adaptive mesh refinement algorithm for finite range interactions based on the a posteriori error estimates \eqref{eq:residual}. Numerical experiments show that our algorithm achieves an optimal convergence rate in terms of accuracy vs. the degrees of freedom, which is also consistent with the optimal a priori error estimates. 

\subsection{Adaptive mesh refinement algorithm.}
\label{sec:numerics:algorithm}

We will develop the adaptive mesh refinement algorithm for GRAC method with finite range interactions. In \cite{APEst:2017}, we have designed the adaptive algorithm for GRAC with nearest neighbour interaction. In the nearest neighbor case, there exists a special set of reconstruction parameters \cite{PRE-ac.2dcorners}. However, with finite range interactions, we have to apply $\ell^{1}$-minimisation to solve the reconstruction parameters from \eqref{eq:E_pt_grac} and \eqref{eq:F_pt_grac}. Furthermore, we need to stabilize the coupling scheme by adding a stabilization term at the interface region, which will effectively add a second order term to the interface stress tensor. We will compare the effects of $\ell^1$ minimization and stabilization in numerical experiments.

\subsubsection{Stress tensor correction}
\label{sec:Interfacial-stress}

By \eqref{eq:defn_Sa}-\eqref{eq:defn_Sh}, the error estimators $\eta_T$, $\eta_M$, and $\eta_C$ depend on the stress tensors $\sh$ and $\sa$, which are unique up to divergence free tensor fields by Lemma \ref{lem:divfree}. In principle, we need to minimize $\eta(y_h) := \tilde{\eta}(\sa(\Ia y_h), \sh(y_h)) = \eta_T(y_h)+\eta_M(y_h)+\eta_C(y_h)$ with respect to all the admissible stress tensors. 

\begin{equation}
		 \< \del\E(\Ia y_h), v \> \leq  \min_{c_a\in N_1(\T_a)^2, c_h\in N_1(\T_h)^2}\tilde{\eta}(\sa(\Ia y_h)+\nabla c_a\mJ, \sh(y_h)+\nabla c_h\mJ)\|\nabla v\|_{L^2}.
		\label{eq:sharpresidual2}
\end{equation}
To save computational power, we can choose a "good'' a/c stress tensor instead of the "optimal" one. A "good" a/c stress tensor should satisfy the following natural conditions:
\begin{itemize}
\item Equal to the atomistic stress tensor in the atomistic domain.
\item Equal to the continuum stress tensor for uniform deformation.
\end{itemize}

Then, we choose $c_a\equiv 0$ and $c_h(q_f) =0$ in \eqref{eq:sharpresidual2}, where $q_f$ is the midpoint of $f\in \Fh, f\cap\L_i=\emptyset$. Also, since $\eta_T$ and $\eta_C$ are higher order contributions to the estimator, we only need to minimize the modeling error $\eta_M$ with respect to the degrees of freedom of $\sh$ adjacent to the interface. 

We propose the following algorithm Algorithm \ref{alg:etc} for approximate stress tensor correction: 


\begin{algorithm}
\label{alg:etc}
Approximate stress tensor correction
\begin{enumerate}
\item Take $\sa(\Ia y_h)$ and $\sh(y_h)$ as the canonical forms in \eqref{eq:defn_Sa} and \eqref{eq:defn_Sh} respectively. 
\item Denote $q_f$ as the midpoint of $f\in \Fh, f\cap\L_i\neq\emptyset$. $c_h$ minimizes the following sum
\begin{equation}
	\sum_{T\in\T^\i}|T|\left[\sa(\Ia y_h, T) - \big(\sh(\Ia y_h, T)+\nabla c_h \mJ\big)\right]^2
\end{equation} 
subject to the constraint that $c_h(q_f)=0$, for $f\bigcap \L_\i = \emptyset$.
\item Let $\sh(y_h)  = \sh(y_h) + \nabla c^h\mJ$, compute $\eta_M$, $\eta_T$ and $\eta_C$ with $\sa(\Ia y_h)$ and $\sh(y_h)$.
\end{enumerate}
\end{algorithm}

\subsubsection{Local error estimator}

We can evaluate the global estimator $\eta$ by computing its three components individually. These components corresond to different operations in the adaptive algorithm. The value of truncation error $\eta_T$ is determined by the domain size. As the size of the computational domain becomes larger, we have smaller $\eta_T$. Hence, we regard $\eta_T$ as a criterion to control the domain size. Values of the modelling error $\eta_M$ and the coarsening error $\eta_C$ indicate local error contributions, and we need to assign $\eta_M$ and $\eta_C$ to local elements properly.  

We define
\begin{displaymath}
\eta_{M}(T, T^{a}):= |T^{a}\bigcap T|\left[\sa(\Ia y_h, T^{a}) - \sum_{T'\in \Th, T'\bigcap T^{a}\neq \emptyset}\frac{|T'\bigcap T^{a}|}{|T^{a}|}(\sh(y_h, T'))\right]^2. 
\end{displaymath}

for $T^{a}\in\Ta$, then let 
\begin{displaymath}
  \eta_{M}(T) = \sum_{T^a\in \Ta, T^a\bigcap T\neq \emptyset} \eta_M(T, T^a), \textrm{ for } T\in\Th. 
\end{displaymath}

\begin{displaymath}
  \eta_{C}(T) = \sqrt{3}C^{\Tr}C'_{\Th}\sum_{f\in \Fh\bigcap T\in\Th}\frac12(h_f\sjump_f)^2.
\end{displaymath}

Once all the local estimators are assigned, we are ready to define the indicator $\rho_{T}$:
\begin{equation}
\rho_{T} = (C^{\Tr})^2\frac{\eta_{M}(T)}{\eta_M}+(\sqrt{3}C^{\Tr}C'_{\Th})^2\frac{\eta_{C}(T)}{\eta_C}.
\label{eq:localestimator1}
\end{equation}

Notice that the sum of local estimators $\rho_T$ together with truncation error $\eta_T$ is equal to the global estimator $\eta$.
The constants $C^{\Tr}$, $C'_{\Th}$ in \eqref{eq:localestimator1} are not known a priorily, instead, we use their empirical estimates in the numerical implementation.

Here is the main algorithm Algorithm \ref{alg:size}, where the D\"{o}rfler adaptive strategy \cite{Dorfler:1996} is adopted.

\def\M{\mathcal{M}}

\begin{algorithm}
\label{alg:size}
A posteriori mesh refinement with size control.
\begin{enumerate}
	\item[Step 0] Set $\Omega_{R_0}$, $\Th$, $N_{\max}$, $\rho_{\textrm{tol}}$, $\tau_1$, $\tau_{3}$ and $R_{\max}$.
	\item[Step 1] \textit{Solve:} Solve the a/c solution $y_{h, R}$ of \eqref{eq:min_ac} on the current mesh $\T_{h, R}$.  
	\item[Step 2] \textit{Estimate}: Carry out the stress tensor correction step in Algorithm \ref{alg:etc}, and compute the error indicator $\rho_{T}$ for each $T\in\Th$, including the contribution from truncation error $\eta_T$. Set $\rho_{T} = 0$ for $T\in\Ta\bigcap \Th$. Compute the degrees of freedom $N$, error estimator $\rho_T$ and $\rho = \sum_{T}\rho_T$. Stop if $N>N_{\max}$ or $\rho < \rho_{\textrm{tol}}$ or $R>R_{\max}$.
	\item[Step 3] \textit{Mark:} 
	\begin{enumerate}
	\item[Step 3.1]:  Choose a minimal subset $\M\subset \Th$ such that
	\begin{displaymath}
		\sum_{T\in\M}\rho_{T}\geq\frac{1}{2}\sum_{T\in\Th}\rho_{T}.
	\end{displaymath}	 
	\item[Step 3.2]: We can find the interface elements which are within $k$ layers of atomistic distance, $\M^k_\i:=\{T\in\M\bigcap \Th^\c: \textrm{dist}(T, \Li)\leq k \}$. Choose $K\geq 1$, find the first $k\leq K$ such that 
	\begin{equation}
		\sum_{T\in \M^k_\i}\rho_{T}\geq \tau_1\sum_{T\in\M}\rho_{T},
		\label{eq:interface2}
	\end{equation}
	
	with tolerance $0<\tau_1<1$.	If such a $k$ can be found, let $\M = \M\setminus \M^k_\i$. Then in step 3, expand interface $\L_\i$ outward by $k$ layers. 
	\end{enumerate}
	\item[Step 4] \textit{Refine:} If \eqref{eq:interface2} is true, expand interface $\L_\i$ outward by one layer. If $\eta_{T}\geq \tau_3 \rho$, enlarge the computational domain. Bisect all elements $T\in \M$. Stop if $\frac{\eta_T}{\eta_M + \eta_C}\geq \tau_2$, otherwise, go to Step 1. 
\end{enumerate}
\end{algorithm}

\subsection{Numerical Experiments}
\label{sec:numerics:problem}
We present the numerical experiments with the following model problem. We consider the two dimensional triangular lattice $\Lhom :=\mA\Z^{2}$ with 
\begin{equation}
\mA = \mymat{1 & \cos(\pi/3) \\ 0 & \sin(\pi/3)}.
\label{eq:Atrilattice}
\end{equation}

Let $a_1 = (1,0)^T$, then $a_j = \mA_6^{j-1}a_1$, $j=1,\dots, 6$, are the nearest neighbour directions in $\Lhom$, where $\mA_6$ is the rotation matrix corresponding to a $\pi/3$ clockwise planar rotation. Let $\Rgnn = \{a_i\}_{i=1}^{6}$ denote the set of nearest neighbor interacting bonds. Given the cut-off radius $\rcut$, then each interaction direction $\rho\in\Rg$ can be uniquely represented as $\rho = \alpha a_{i} + \beta a_{i+1}$, where $\alpha\geq 0, \beta\geq 0, \alpha+\beta\leq\rcut$ and $a_i\in \Rgnn$. 

Recall the EAM potential defined in \eqref{eq:eam_potential}. Let
\begin{displaymath}
\phi(r)=\exp(-2a(r-1))-2\exp(-a(r-1)),\quad \psi(r)=\exp(-br)
\end{displaymath}
\begin{displaymath}
F(\tilde{\rho})=C\left[(\tilde{\rho}-\tilde{\rho}_{0})^{2}+
(\tilde{\rho}-\tilde{\rho}_{0})^{4}\right]
\end{displaymath}
with parameters $a=4.4, b=3, c=5$ and $\tilde{\rho}_{0}=6\exp(-b)$, which is the same as the numerical experiments in the a priori analysis of GRAC method \cite{COLZ2013}.

To generate a defect, we remove $k$ atoms from $\Lhom$,
\begin{align*}
\L_{k}^{\textrm{def}}:=\{-(k/2)e_{1}, \ldots, (k/2-1)e_{1})\},     & \qquad{\textrm{if}}\quad k \quad\textrm{is even},\\
\L_{k}^{\textrm{def}}:=\{-(k-1)/2e_{1}, \ldots, (k-1)/2e_{1})\}, & \qquad{\textrm{if}}\quad k \quad\textrm{is odd},
\end{align*}
and $\L = \Lhom\setminus \L_{k}^{\textrm{def}}$. 

For $\ell\in\L$, consider the next nearest neighbour interaction, $\Nhd_{\ell} := \{ \ell' \in \L \sep 0<|\ell'-\ell| \leq 2 \}$, and interaction range $\Rg_\ell := \{\ell'-\ell \sep \ell'\in \Nhd_\ell\} \subseteq \{a_j, j=1,\dots, 18\}$. The defect core $\Ddef$ can be defined by $\Ddef = \{x: \textrm{dist} (x,  \L_{k}^{\textrm{def}})\leq 2 \}$, $\L \bigcap \Ddef$ is the first layer of atoms around $\L_{k}^{\textrm{def}}$.

\subsubsection{Di-vacancy}
\label{sec:numerics:di-vacancy}

In this section, we numerically justify the performance of the proposed adaptive mesh refinement algorithm.
We take the same di-vacancy example in \cite{COLZ2013}, namely, setting $k=2$ for $\L_{k}^{\textrm{def}}$. We apply isotropic stretch $\mathrm{S}$ and shear $\gamma_{II}$  by setting
\begin{displaymath}
{\textsf{B}}=\left(
	\begin{array}{cc}
		1+\mathrm{S} & \gamma_{II} \\
		0            & 1+\mathrm{S}
	\end{array}	 \right)
	\cdot{\mathsf{F_{0}}}
\end{displaymath}
where $\mathsf{F_{0}} \propto \mathrm{I}$ minimizing the Cauchy-Born energy density $\mathrm{W}$, $\mathrm{S}=\gamma_{II}=0.03$.

\subsubsection{Micro-crack}
\label{sec:micro-crack}

In the microcrack experiment, we remove a longer segment of atoms,
$\L^{\textrm{def}}_{11}=\{-5e_1,\dots,5e_1\}$ from the computational
domain. The body is then loaded in mixed mode ${\rm I}$ \& ${\rm II}$,
by setting,
\begin{displaymath}
\mB := \begin{pmatrix} 1 & \gamma_{\rm II} \\ 0 & 1+\gamma_{\rm I} \end{pmatrix}\cdot \mF_0.
\end{displaymath}
where $\mF_0 \propto I$ minimizes $W$, and $\gamma_{\rm I}=\gamma_{\rm
  II}=0.03$ ($3\%$ shear and $3\%$ tensile stretch).  

We take $\tau_3 = 0.7$ in the numerical implementation of Algorithm \ref{alg:size}. We compare results with and without stabilisation and also with different optimisation approaches to obtain the reconstruction parameters. 

\begin{figure}[H]
\begin{center}
	\includegraphics[scale=0.55]{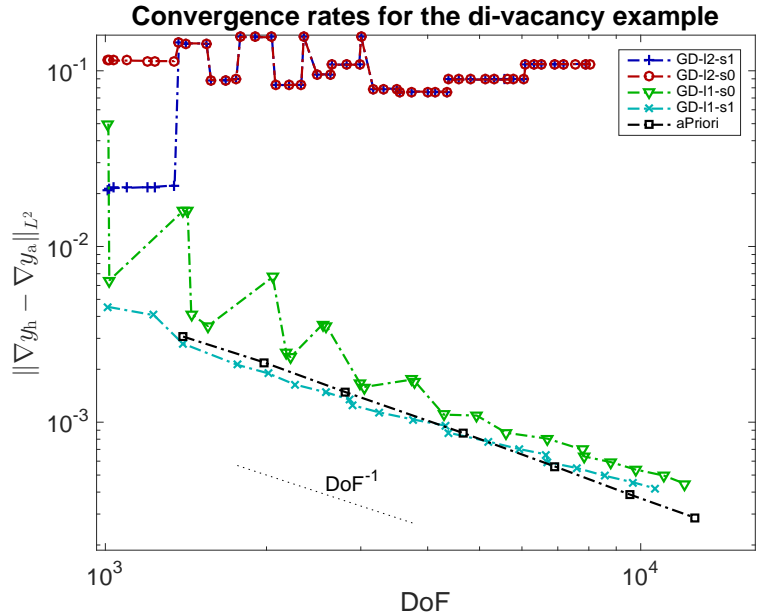}
	\caption{Numerical results by Algorithm \ref{alg:size}: $H^1$ error vs. $N$, using a posteriori estimator in $H^1$ norm. In the legend, l1 means $\ell^{1}$-minimization approach while l2 represents a least-squares approach; s1 indicates with stabilization and s0 without.}
	\label{fig:nV2-H1}
\end{center}
\end{figure}

\begin{figure}[H]
\begin{center}
	\includegraphics[scale=0.55]{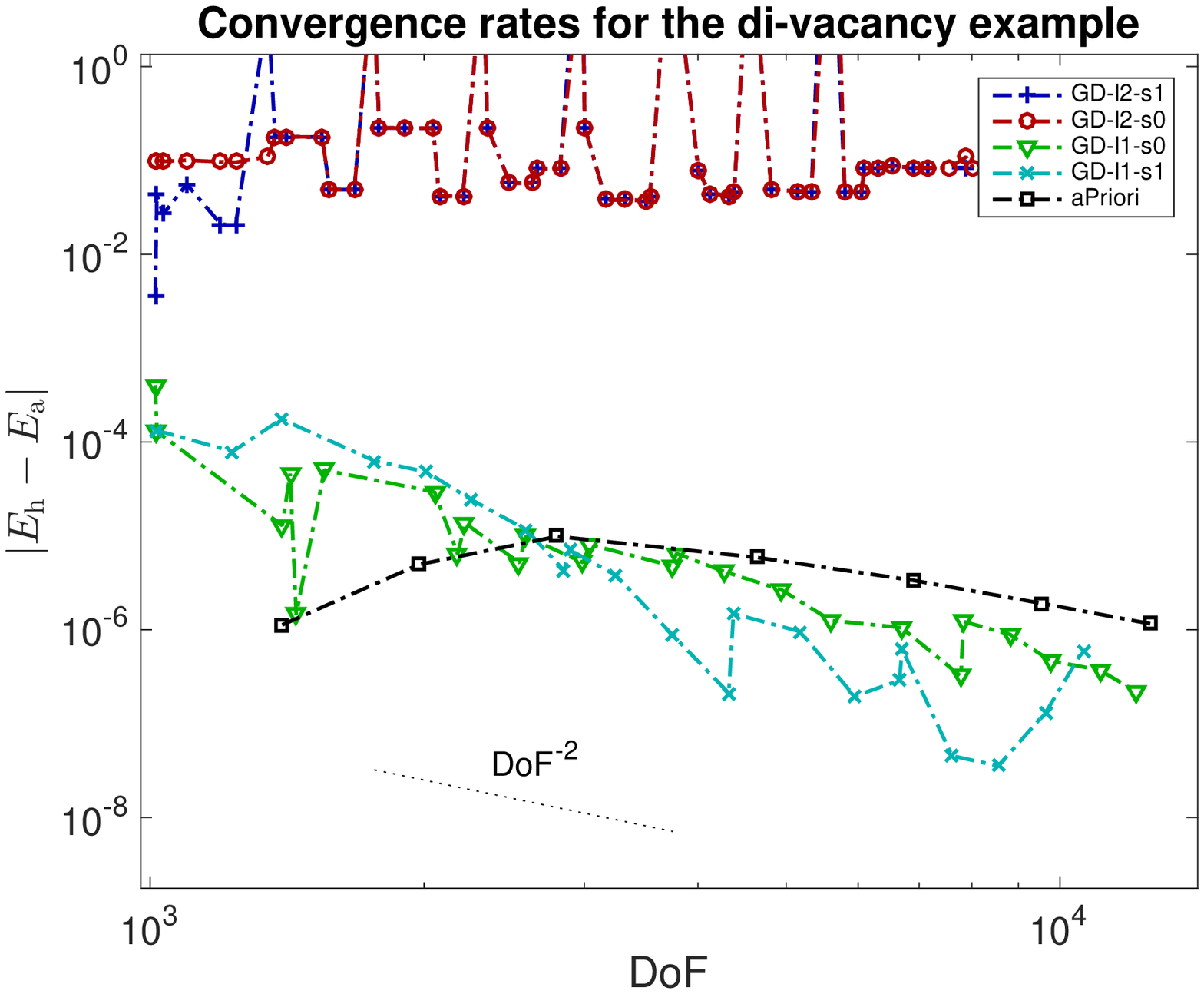}
	\caption{Numerical results by Algorithm \ref{alg:size}: Energy error vs. $N$, using a posteriori estimator in $H^1$ norm. In the legend, l1 means $\ell^{1}$-minimization approach while l2 represents a least-squares approach; s1 indicates with stabilization and s0 without.}
	\label{fig:nV2-En}
\end{center}
\end{figure}

\begin{figure}[H]
\begin{center}
	\includegraphics[scale=0.55]{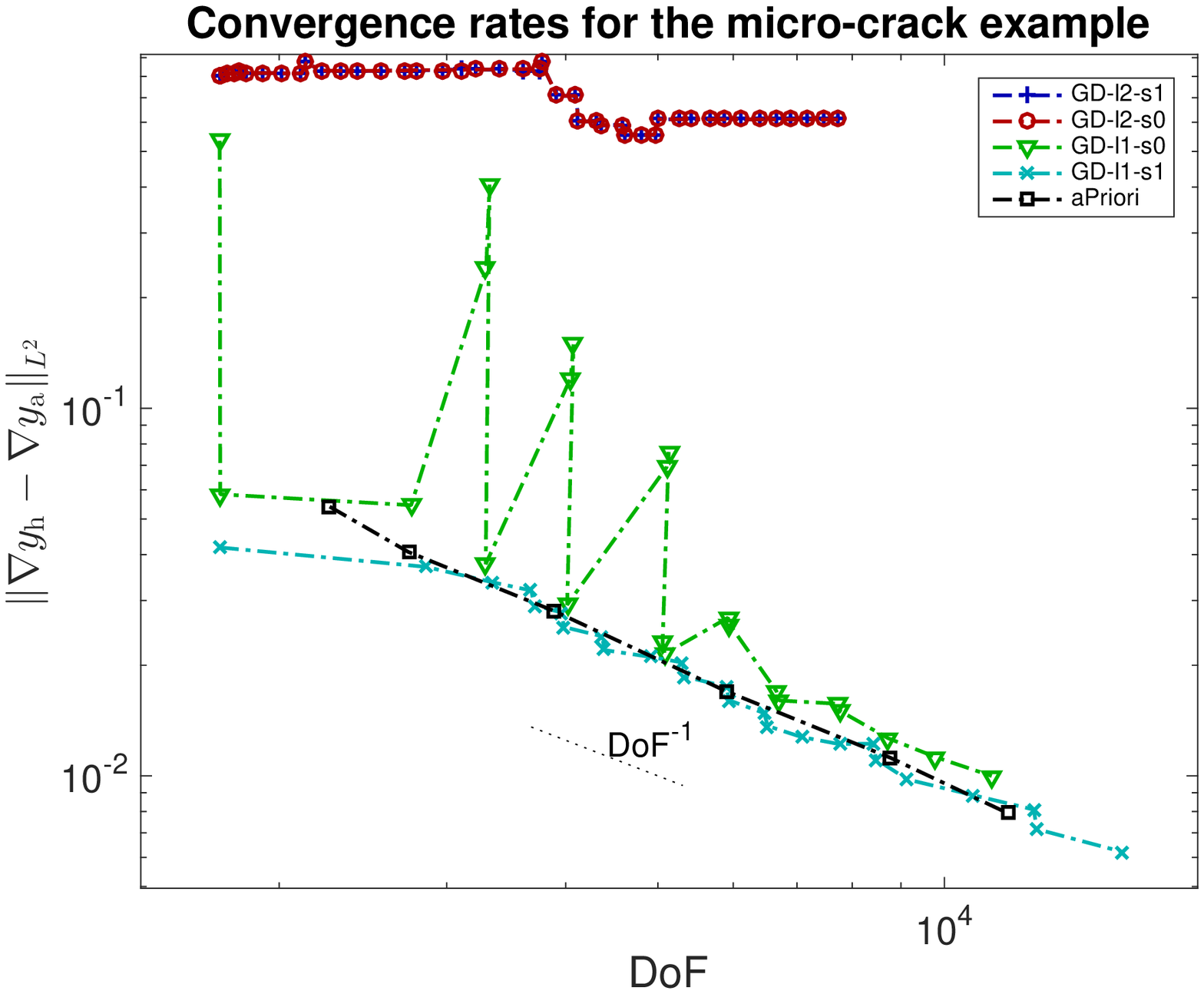}
	\caption{Numerical results by Algorithm \ref{alg:size}: $H^1$ error vs. $N$, using a posteriori estimator in $H^1$ norm. In the legend, l1 means $\ell^{1}$-minimization approach while l2 represents a least-squares approach; s1 indicates with stabilization and s0 without.}
	\label{fig:nV11-H1}
\end{center}
\end{figure}

\begin{figure}[H]
\begin{center}
	\includegraphics[scale=0.55]{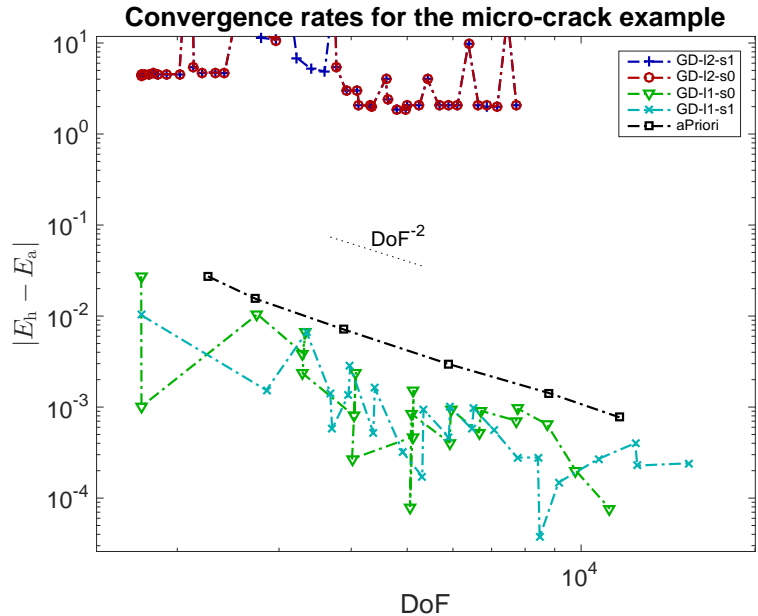}
	\caption{Numerical results by Algorithm \ref{alg:size}: Energy error vs. $N$, using a posteriori estimator in $H^1$ norm. In the legend, l1 means $\ell^{1}$-minimization approach while l2 represents a least-squares approach; s1 indicates with stabilization and s0 without.}
	\label{fig:nV11-En}
\end{center}
\end{figure}

From the numerical results in Figures \ref{fig:nV2-H1} - \ref{fig:nV11-En}, we can see that the least square method for the reconstruction parameters does not converge at all. For the $\ell^1$-minimisation approach without stabilization, there exist large errors in the pre-asymptotic regime, though the solutions tend to converge with increasing degree of freedom. For $\ell^1$-minimisation combined with stabilisation, we are able to obtain the optimal convergence rate that consistent with the a priori result. 

\section{Conclusion}
\label{sec:conclusion}
In this paper, we construct the adaptive algorithm for a class of consistent (ghost force free) atomistic/continuum coupling schemes with finite range interactions based on the rigorous a posteriori error estimates. Different from the localization formula for the stress tensor in the a priori analysis \cite{Or:2011a, OrShSu:2012, OrtnerTheil2012}, we develop a computable formulation for the stress tensor. Combined with $\ell_1$-minimization approach for reconstruction parameters and stabilization, we have shown that the numerical results for the corresponding adaptive algorithms are comparable to optimal a priori analysis. 

The extension of the straight screw dislocation in two dimensions and point defect case in three dimensions is straightforward. More practical problems, for example, the study of dislocation nucleation and dislocation interaction by a/c coupling methods has attracted considerable attention from the early stage of a/c coupling methods \cite{Tadmor:1999, Phillips:1999}. The difficulty is to deal with boundary condition and complicated geometry changes of the interface.

For general atomistic/continuum coupling schemes, such as BQCE, BQCF and BGFC, the a priori analysis in \cite{MiLu:2011, LiOrShVK:2014, OrZh:2016} provide a general analytical framework and the stress tensor based formulation plays a key role in the analysis. Therefore, the a posteriori analysis for those coupling schemes can inherit this analytical framework and the stress tensor formulation. Techniques developed in this paper will be essential for the efficient implementation of the corresponding adaptive algorithms.


\section{Appendix}
\label{sec:appendix}

\subsection{Derivation of \eqref{eq:stress:drho} for Two Dimensional Triangular Lattice}
\label{sec:stress:assm}

In this section, for two dimensional triangular lattice introduced in \S~\ref{sec:numerics:problem}, we derive \eqref{eq:stress:drho}.
We first classify interaction bonds into two types according to their intersection to corresponding elements. 


\def\Rgnn{\Rg_{\textrm{nn}}}
\def\mBp{\mB_\textrm{I}}
\def\mBc{\mB_\textrm{II}}

\begin{figure}[H]
\begin{center}
	\includegraphics[scale=0.45]{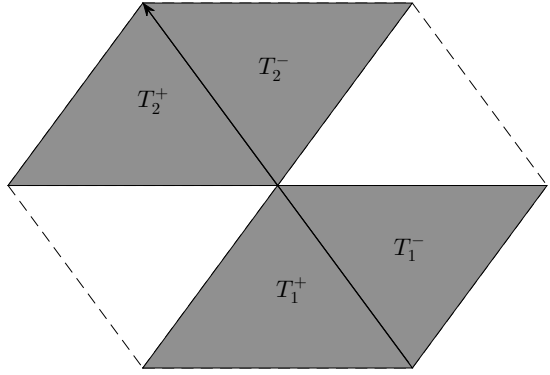}
	\caption{Illustration of corresponding elements (marked gray) to bond $\rho\in\mBp$ with $\rho = 2a_3$.}
	\label{figs:a2b0}
\end{center}
\end{figure}

Type I interaction bonds is parallel to one of the nearest neighbour directions. The set of all Type I bonds is defined by $\mBp = \{\rho| \exists \a_i \in \Rgnn, \text{ such that }\rho = |\rho| a_i\}$, as illustrated in Figure \ref{figs:a2b0}.

For each $\rho\in\mBp$ with start point $\ell\in\L$ and $v\in\Use$, we have
\begin{align}
 D_\rho v &= \sum_{k=1}^{|\rho|} D_{a_i}v(\ell+(k-1)a_i)  \nonumber \\
 & =  \left(\sum_{k=1}^{|\rho|} \left(\frac12\D_{T_{k}^{+}}v+\frac12\D_{T_{k}^{-}}v\right)\right)\cdot \frac{\rho}{|\rho|} \\ 
 & =  \sum_{k=1}^{|\rho|} \frac{1}{2|\rho|} \left(\D_{T_{k}^{+}}v\cdot \rho+\D_{T_{k}^{-}}v\cdot \rho\right) \nonumber
\end{align}
where $T_{k}^{+}$ and $T_{k}^{-}$ are the triangles with $\left(\ell+(k-1)a_i, \ell+k a_i\right)$ as the common edge,  for $k=1,\ldots,|\rho|$, superscript ``$+$'' means the element being located on the left of the bond. In this case, the contribution in\eqref{eq:stress:drho} is 
$$\wlrho(T_{k}^{(\pm)})=\frac{1}{2|\rho|}.$$

\begin{figure}[H]
\begin{center}
	\includegraphics[scale=0.35]{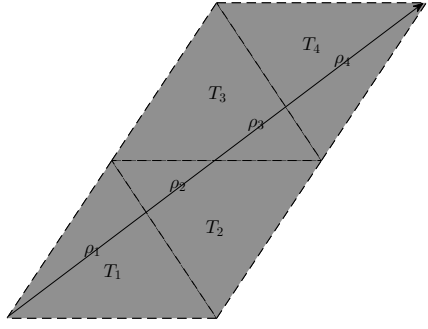}
	\caption{Illustration of the intersecting elements (marked gray) with respect to bond $\rho\in\mBc$, $\rho = a_1 + 2a_2$. }
	\label{figs:a1b2}
\end{center}
\end{figure}

Type II interaction bonds is not parallel to any of the nearest neighbour directions, therefore they must cross the intersecting elements. The set of all Type II bonds is defined by $\mBc = \{\rho|\rho \text{ is not parallel to any }a_i\in\Rgnn\}$, as illustrated in Figure \ref{figs:a1b2}.

For any Type II interaction bond $\rho\in\mBc$, there exists nearest neighbour directions $a_i$, $a_{i+1}$, such that $\rho = \alpha a_i + \beta a_{i+1}$,  $\alpha,\beta>0$. Let $(\ell, \ell+\rho)$ consecutively intersect with elements $T_k$, $k = 1, \cdots, n_p$, for $v\in\Use$, we have
\begin{align}
 D_\rho v &= \partial_\rho V_\ell \left(\sum_{k=1}^{n_\rho} D_{\rho_{k}}v(\ell+\sum_{t=1}^{k}\rho_{k-1})\right)  \nonumber \\
 & =  \left(\sum_{k=1}^{n_\rho} \D_{T_k}v\right)\cdot \rho_k \\ 
 & =  \left(\sum_{k=1}^{n_\rho} \D_{T_k}v\right)\cdot \frac{|\rho_k|}{|\rho|}\rho \nonumber \\ 
 & =  \sum_{k=1}^{n_\rho)} \frac{|\rho_k|}{|\rho|} \D_{T_k}v\cdot \rho  \nonumber
\end{align}
where the set of vectors $\{\rho_k\}_{k=1}^{n_\rho}$ is the unit partition of $\rho$ with $\rho = \sum_{k=1}^{n_\rho}\rho_k$, $\rho_0 = \mathbf{0}$ and $\frac{\rho_k}{|\rho_k|}=\frac{\rho}{|\rho|}$ for $k=1,\ldots,n_p$. $\rho_k = T_k \cap (\ell, \ell+\rho)$. For the two dimensional triangular lattice as in \S~\ref{sec:numerics:problem}, we have $n_\rho = 2(\alpha + \beta -1)$ if $\alpha \neq \beta$, and $n_\rho = 2\alpha$ if $\alpha = \beta$. In this case, 
$$\wlrho(T_k)=\frac{|\rho_k|}{|\rho|}.$$

\begin{figure}[H]
\begin{center}
	\includegraphics[scale=0.45]{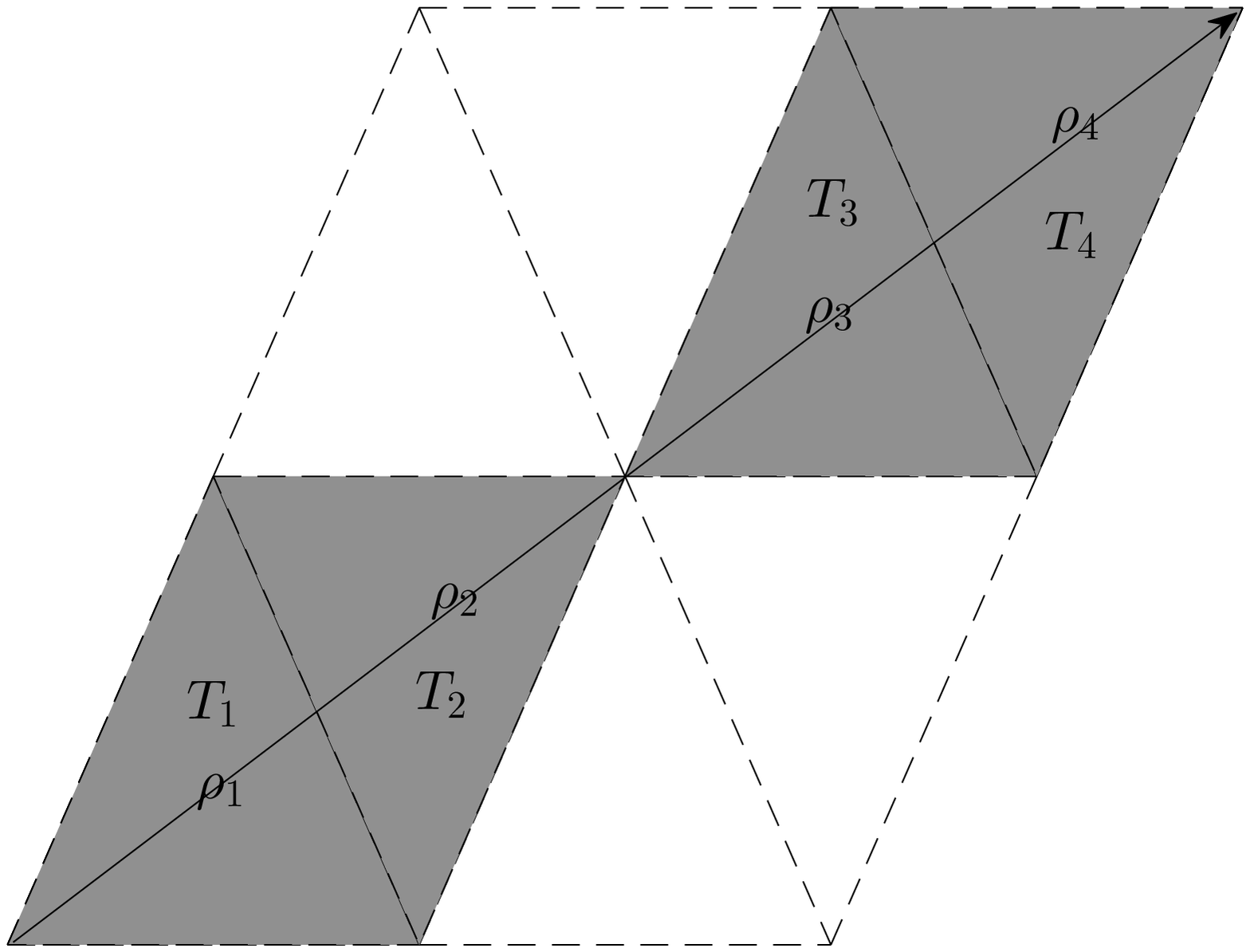}
	\caption{Illustration of corresponding elements (marked gray) to bond $\rho\in\mBp$ with $\rho = 2a_1+2a_2$.}
	\label{figs:a2b2}
\end{center}
\end{figure}

For the Type II bonds with $\alpha = \beta$, in the hexagon site geometry, The contribution factor $\omega_\rho^{T}$ is the same for every element related to bond $\rho$ with the value $\omega_\rho^T = \frac{1}{2\alpha}$ which can be see from figure (\ref{figs:a2b2}).

%


\bibliographystyle{plain}
\bibliography{qc-1.bib}

\end{document}